\documentclass[letterpaper, paper,11pt]{AAS}		

\usepackage{bm}
\usepackage{amsmath}
\usepackage{amssymb}
\usepackage{amsthm}

\theoremstyle{remark}
\newtheorem{remark}{Remark}
\usepackage{graphicx}
\usepackage[colorlinks=true, pdfstartview=FitV, linkcolor=black, citecolor= black, urlcolor= black]{hyperref}
\usepackage{overcite}
\usepackage{footnpag}	
\usepackage{algorithm}
\usepackage{algpseudocode}
\usepackage{subcaption}

\PaperNumber{26-941}

\begin{document}

\title{An Efficient iLQR Algorithm for Optimal Control with Bounded States and Inputs: Application to Spacecraft Maneuvers}

\author{Abhijeet\thanks{Graduate Student, Aerospace Engineering, Texas A\& M University, College Station, TX.},  
Tarun Hejmadi\thanks{Graduate Student, Aerospace Engineering, Texas A\& M University, College Station, TX.},
\ and Suman Chakravorty\thanks{Professor, Aerospace Engineering, Texas A\& M University, College Station, TX.}
}

\maketitle{}

\begin{abstract}
This paper presents a constrained iterative Linear Quadratic Regulator (Box-iLQR) method for nonlinear aerospace optimal-control problems with bounded states and controls. Logarithmic barrier functions are used with a two-fold objective: the barrier parameter is reduced toward zero to recover the constrained optimal trajectory, while a finite barrier is retained to construct a smooth, constraint-aware feedback policy for closed-loop operation. The method is demonstrated on spacecraft attitude control, a minimum-fuel planar orbit transfer, and an \(L_2\)-to-\(L_2\) halo-orbit transfer. Box-iLQR respects the imposed constraints and produces nominal trajectories together with time-varying feedback policies. Closed-loop simulations under disturbance, measurement noise, and control uncertainty show substantial reductions in terminal error. These results demonstrate that Box-iLQR provides both high-accuracy constrained trajectory optimization and practical constraint-aware feedback within a single framework.

\end{abstract}

\section{Introduction}
\label{sec:introduction}

Constrained trajectory optimization is fundamental to a broad range of
aerospace guidance, navigation, and control problems. Representative
applications include spacecraft attitude reorientation, rendezvous and
proximity operations, low-thrust orbit transfer, planetary landing, launch
vehicle guidance, and atmospheric entry. In these applications, the control
trajectory must optimize a mission-dependent performance index while
simultaneously satisfying nonlinear dynamics, prescribed boundary conditions,
and path constraints arising from actuator saturation, pointing restrictions,
collision avoidance, structural limits, and operational requirements. Because
the governing dynamics are generally nonlinear and the constraints may become
active over extended portions of the trajectory, analytical solutions are
available only for a limited number of special cases. Consequently, reliable
and computationally efficient numerical optimal-control methods are essential
for both mission design and onboard trajectory generation
\cite{betts1998survey,rao2009survey,kelly2017introduction,malyuta2021advances}.

Numerical optimal-control methods are commonly divided into indirect and direct
approaches. Indirect methods apply the calculus of variations or Pontryagin's
minimum principle to derive the necessary conditions for optimality. The
resulting state--costate boundary-value problem can provide highly accurate
solutions and direct information about the optimal switching structure.
However, the method may be sensitive to the initialization of the unknown
costates, particularly for long-duration trajectories, strongly nonlinear
dynamics, and problems containing discontinuous or nearly discontinuous
controls. Direct methods instead discretize the state and control trajectories
and transcribe the optimal-control problem into a finite-dimensional nonlinear
program. These methods are generally easier to initialize and can accommodate
a broad class of path and boundary constraints, but the resulting nonlinear
program may become large as the number of discretization points increases
\cite{betts1998survey,rao2009survey,kelly2017introduction}. These considerations
have motivated considerable interest in trajectory-optimization algorithms
that retain the computational structure of dynamic programming while avoiding
the solution of a large, general-purpose nonlinear program \cite{Lantoine2012HybridDDP,Lantoine2012HybridDDPApplication}.

Spacecraft attitude maneuvering provides a representative example of a
nonlinear constrained optimal-control problem. Large-angle reorientation
maneuvers involve nonlinear coupling between the attitude kinematics and
rigid-body rotational dynamics, while practical actuators impose strict bounds
on the available torque and angular momentum. The classical work of Vadali and
Junkins considered both optimal open-loop maneuvers and stabilizing feedback
control for rigid spacecraft
\cite{vadali1984optimalopenloop}. Nonlinear optimal-feedback approximations for
large-angle maneuvers were subsequently developed by Carrington and Junkins
\cite{carrington1986optimal}. More recently, discrete-mechanics-based
optimal-control formulations have been proposed to explicitly account for
momentum and control constraints while preserving the underlying attitude
geometry \cite{phogat2018discretetime}. These studies demonstrate the importance
of incorporating actuator bounds into the maneuver-design process rather than
applying saturation to an unconstrained control solution after optimization.

Fuel-optimal low-thrust orbit transfer presents a related but distinct
computational challenge. The long transfer duration, oscillatory orbital
dynamics, terminal constraints, and limited thrust authority can produce a
highly sensitive optimization problem. Moreover, minimum-fuel formulations
typically generate bang--off--bang thrust histories, for which the number and
location of switching times are not known a priori. Homotopy methods have been
used to continuously transform a smoother optimal-control problem into the
desired minimum-fuel problem \cite{haberkorn2004lowthrust}, while
pseudospectral discretizations have been employed to obtain high-accuracy
low-thrust solutions \cite{ross2007lowthrust}. Taheri, Kolmanovsky, and Atkins
developed a smoothing procedure for indirect minimum-fuel trajectory
optimization that regularizes the discontinuous thrust law and progressively
recovers the bang--off--bang solution
\cite{taheri2016enhancedsmoothing}. Although these methods can produce accurate
solutions, the treatment of nonsmooth controls, active control bounds, and
robust closed-loop execution remains an important consideration for autonomous
applications.

Differential dynamic programming and the iterative Linear Quadratic Regulator
(iLQR) provide an alternative class of numerical optimal-control methods based
on successive local approximations of the nonlinear problem. The underlying second-order trajectory-improvement principle was established
in the early work of Mayne \cite{mayne1966secondorder}. Recent work has shown
that, away from a solution, omitting the second-order dynamics terms can provide
a more consistent local approximation, motivating the use of iLQR rather than
full DDP in the present work
\cite{abhi_DDP_ILQR,abhi_DDP_ILQR_ACC}. The iLQR formulation uses a local
linear approximation of the dynamics and a quadratic approximation of the
cost to compute an improved control policy through alternating backward and
forward passes \cite{li2004iterative}. Because each forward rollout satisfies
the nonlinear dynamics, iLQR is naturally suited to shooting-based trajectory
optimization. Furthermore, the backward pass produces both a feedforward
control correction and a time-varying state-feedback gain. The resulting
policy can therefore be used not only to improve the nominal trajectory but
also to stabilize its execution. The ability of this approach to generate and
stabilize complex nonlinear motions has been demonstrated in online
trajectory-optimization applications
\cite{tassa2012synthesis}.

Standard iLQR is formulated for unconstrained state and control spaces.
Consequently, directly applying the algorithm to aerospace systems can produce
controls that exceed actuator capabilities or trajectories that violate
mission-critical state restrictions. Control-limited differential dynamic
programming addresses box constraints on the control by solving a local
bound-constrained quadratic program during the backward pass
\cite{tassa2014controllimited}. However, the simultaneous treatment of state
and control constraints is more difficult because state feasibility depends
on the complete nonlinear forward propagation. Simple quadratic penalties do
not guarantee feasibility and may require large penalty coefficients that
degrade the conditioning of the local quadratic model. Methods such as ALTRO, FATROP, and primal--dual interior-point solvers such as IPOPT primarily target the computation of constrained optimal trajectories; although Riccati or iLQR recursions may be used internally, the resulting constrained solution does not inherently provide a constraint-aware feedback policy suitable for closed-loop execution near active constraints \cite{Howell2019ALTRO,Vanroye2023FATROP,Wachter2006IPOPT}.

In this work, state and control inequality constraints are incorporated into
iLQR through logarithmic barrier functions. Interior-point and barrier methods
are well-established tools in numerical optimization for treating inequality
constraints by replacing the original constrained problem with a sequence of
smooth unconstrained subproblems
\cite{nocedal2006numerical}. For the present formulation, the barrier terms are
included directly in the running and terminal costs. In addition to
discouraging the trajectory from approaching an inadmissible region, the
second derivatives of the logarithmic barriers contribute positive curvature
to the local quadratic approximation used in the iLQR backward pass. This
curvature acts as a natural regularization mechanism and is particularly
useful for minimum-fuel problems in which the running cost may have little or
no curvature over coast intervals
\cite{abhijeet2026safeoptimalcontrolusing}.

The resulting method, referred to herein as Box-iLQR, uses a nested solution
procedure. For fixed barrier parameters, an inner iLQR loop performs backward
and forward passes with regularization and line search until the corresponding
barrier subproblem converges. An outer continuation loop then reduces the
barrier parameters so that the sequence of strictly feasible solutions
approaches a solution of the original constrained problem. The solution of
each barrier subproblem is used to warm-start the next one. When a barrier
reduction is too aggressive and causes the inner solve to fail, the reduction
factor associated with the corresponding constraint is adjusted to obtain a
more conservative continuation step. This adaptive relaxation strategy avoids
requiring a single barrier schedule to be selected a priori for constraints
with different scales and sensitivities.

An important feature of the proposed approach is that the barrier parameter
plays two distinct roles. Barrier continuation toward zero is used to approach
the constrained nominal optimum, while the feedback gains obtained at a
selected finite barrier value are retained for closed-loop execution. The
finite barrier preserves constraint-dependent curvature in the local quadratic
model and provides a locally stabilizing, constraint-aware feedback policy.
This is particularly useful for spacecraft applications subject to actuator
uncertainty, environmental disturbances, navigation errors, and modeling
uncertainty.

The principal contributions of this work are summarized as follows:
\begin{enumerate}
    \item A logarithmic-barrier iLQR formulation is presented for nonlinear
    optimal-control problems subject to componentwise state and control
    bounds.

    \item An adaptive barrier-relaxation procedure is introduced to reduce
    sensitivity to the prescribed continuation schedule while maintaining a
    feasible warm start for successive barrier subproblems.

   \item The time-varying feedback gains obtained at a selected finite barrier
    value are retained for closed-loop execution, providing a constraint-aware
    feedback policy without requiring a separate controller-design procedure.

   \item The method is evaluated on three spacecraft problems: a fully actuated
    large-angle attitude maneuver, a minimum-fuel planar orbit transfer, and a
    fixed-time $L_2$-to-$L_2$ Halo-orbit transfer in the Earth--Moon CR3BP.
    Closed-loop simulations assess robustness to disturbance torques, measurement
    noise, and control uncertainty.
\end{enumerate}

The remainder of this paper is organized as follows. First, the continuous-
and discrete-time constrained optimal-control problems are defined. The
logarithmic barrier formulation and the Box-iLQR backward and forward passes
are then developed, followed by the adaptive barrier-relaxation algorithm.
Numerical results are subsequently presented for spacecraft attitude control, a minimum-fuel planar orbit transfer, and a fixed-time $L_2$-to-$L_2$ Halo-orbit transfer in the Earth--Moon CR3BP, including closed-loop robustness studies.

\section{Background: The Optimal Control Problem}
Consider the following optimal control problem:
\begin{subequations} \label{eq:ocp}
    \begin{align}
        \min_{\mathbf{u}(\cdot)} \quad & \mathcal{J} = \phi(\mathbf{x}(t_f)) + \int_{0}^{t_f} c(\mathbf{x}(t),\mathbf{u}(t))\,dt  \\
        \text{s.t.} \quad & \dot{\mathbf{x}}(t) = \mathbf{f}(\mathbf{x}(t),\mathbf{u}(t);t) \quad \forall ~t \in [0,t_f],\label{eq:ocp_dyn} \\
        & \mathbf{x}(t) \in \mathcal{X} := \{\mathbf{x} \in \mathbb{R}^n \mid \mathbf{\underline{x}} \leq \mathbf{x} \leq \mathbf{\bar{x}}\} \label{eq:ocp_x_bounds}, \\
        & \mathbf{u}(t) \in \mathcal{U} := \{\mathbf{u} \in \mathbb{R}^m \mid \mathbf{\underline{u}} \leq \mathbf{u} \leq \mathbf{\bar{u}}\} \label{eq:ocp_u_bounds}, \\
        & \mathbf{x}(0) ~~\text{and}~~ t_{f} ~~\text{are given.}\nonumber
    \end{align}
\end{subequations}
In the above equation, $\mathbf{x}(t) \in \mathbb{R}^n$ is the state and $\mathbf{u}(t) \in \mathbb{R}^m$ is the control input. The function $\mathbf{f}:\mathbb{R}^n \times \mathbb{R}^m \times \mathbb{R} \to \mathbb{R}^n$ describes the system dynamics. The objective functional $\mathcal{J}$ is composed of a terminal cost $\phi:\mathbb{R}^n \to \mathbb{R}$ and a running cost $c:\mathbb{R}^n \times \mathbb{R}^m \to \mathbb{R}$. The state and control are constrained to lie within the compact sets $\mathcal{X}$ and $\mathcal{U}$, respectively, for all $t \in [0, t_f]$. 

The discrete counterpart of the problem mentioned in eq. \eqref{eq:ocp} is:
\begin{subequations} \label{eq:ocp_dis}
    \begin{align}
        \min_{\mathbf{u}_t} \quad & \mathcal{J} = \Phi(\mathbf{x}_T) + \sum_{t=0}^{T-1} C(\mathbf{x}_t,\mathbf{u}_t)\ , \\
        \text{s.t.} \quad & {\mathbf{x}}_{t+1} = \mathbf{F}(\mathbf{x}_t,\mathbf{u}_t;t) \quad \forall~ t=0,1\cdots T-1, \label{eq:ocp_dyn_dis} \\
        & \mathbf{x}_t \in \mathcal{X} := \{\mathbf{x} \in \mathbb{R}^n \mid \mathbf{\underline{x}} \leq \mathbf{x} \leq \mathbf{\bar{x}}\} \label{eq:ocp_x_bounds_dis}, \\
        & \mathbf{u}_t \in \mathcal{U} := \{\mathbf{u} \in \mathbb{R}^m \mid \mathbf{\underline{u}} \leq \mathbf{u} \leq \mathbf{\bar{u}}\}, \label{eq:ocp_u_bounds_dis}\\
        & \mathbf{x}_{0} ~,~ T ~~\text{are given} \nonumber
    \end{align}
\end{subequations}
and $\Phi, C(\cdot),$ and $\mathbf{F}(\cdot, \cdot)$ are the discrete equivalents of $\phi, c$ and $\mathbf{f}$, respectively.

We can incorporate the inequality constraints into the cost function using a logarithmic barrier method. The constrained optimal control problem in \eqref{eq:ocp_dis} is thereby transformed into an unconstrained problem by augmenting the objective with barrier terms and the dynamic constraints with Lagrange multipliers. This leads to the formulation of the following barrier subproblem, where we seek the stationary point of the augmented cost function $\mathcal{J}_{\boldsymbol{\mu},\boldsymbol{\sigma}}$:
\begin{align} \label{eq:augmented_cost}
    &\mathcal{J}_{\boldsymbol{\mu},\boldsymbol{\sigma}} = \Phi(\mathbf{x}_T) -  \sum_{i \in \mathcal{I}_x} \mu_{i}\left( \log(x_{T,i} - \underline{x}_{i}) + \log(\bar{x}_i - x_{T,i}) \right) \nonumber \\
    &\quad + \sum_{t=0}^{T-1} \Bigg[ C(\mathbf{x}_t, \mathbf{u}_t) -  \sum_{i \in \mathcal{I}_x} \mu_{i}\left( \log(x_{t,i} - \underline{x}_{i}) + \log(\bar{x}_i - x_{t,i}) \right) \nonumber \\
    &\quad\quad -  \sum_{j \in \mathcal{I}_u} \sigma_{j}\left( \log(u_{t,j} - \underline{u}_{j}) + \log(\bar{u}_j - u_{t,j}) \right) + \boldsymbol{\lambda}_{t+1}^{\top} (\mathbf{F}(\mathbf{x}_t, \mathbf{u}_t) - \mathbf{x}_{t+1}) \Bigg].
\end{align}
In the expression above, the variables are defined as follows:
$\boldsymbol{\mu} = [\mu_{1},\mu_{2},\cdots\mu_{n_1}]$ $~~ \text{s.t.}~~\mu_{i} > 0~~\forall i\in1,\cdots n_{1}$ and $\boldsymbol{\sigma} = [\sigma_{1},\sigma_{2},\cdots\sigma_{m_1}]~~ \text{s.t.}~~\sigma_{j} > 0~~\forall j\in 1,\cdots m_{1}$ are the barrier parameters corresponding to the state and control inequality constraints, respectively.

\begin{remark}
The discrete-time problem in Eq.~\eqref{eq:ocp_dis} provides a numerical transcription of the continuous-time problem in Eq.~\eqref{eq:ocp} for a sufficiently small discretization interval $\Delta t$. For state constraints, feasibility at the discrete nodes does not, in general, guarantee feasibility between nodes. A sufficient condition is that each constrained state component
$x_i(t)$ be monotonic over every interval $[t_k,t_{k+1}]$. In this case, its extrema occur at the interval endpoints, and satisfaction of the corresponding bounds at $t_k$ and $t_{k+1}$ implies satisfaction throughout the interval. When this condition does not hold, intersample constraint satisfaction must be
verified or enforced separately.
\end{remark}

\section{Theory: the Box-iLQR algorithm}
The algorithm is initialized with a feasible trajectory. We note that finding such a trajectory can be a challenging problem in itself, but methods for doing so are outside the scope of this work. The core of the algorithm consists of the backward and forward passes, which we detail below.

Our cost function is augmented with logarithmic barrier terms to enforce inequality constraints. The running cost includes a barrier term $\boldsymbol{\omega}(\mathbf{x}_t, \mathbf{u}_t)$, and the terminal cost includes a barrier term $\boldsymbol{\Omega}(\mathbf{x}_T)$.

The \textbf{running barrier term} is defined as:
\begin{align}
\label{eq:running_barrier}
\boldsymbol{\omega}(\mathbf{x}_t, \mathbf{u}_t) = &- \sum_{i \in \mathcal{I}_x} \mu_{i}\left( \log(x_{t,i} - \underline{x}_{i}) + \log(\overline{x}_i - x_{t,i}) \right) \nonumber\\
&- \sum_{j \in \mathcal{I}_u} \sigma_{j}\left( \log(u_{t,j} - \underline{u}_{j}) + \log(\overline{u}_j - u_{t,j}) \right)
\end{align}
Its first-order derivative (gradient) will have some zero and some non-zero elements. The second-order derivatives (Hessians) is a diagonal matrix, with some zero and some non-zero elements. The non-zero elements of gradient and Hessian are given by:
\begin{align}
\label{eq:log_derivatives}
    &(\boldsymbol{\omega}_{\mathbf{x}})_i = -\mu_i \left( \frac{1}{x_{t,i} - \underline{x}_i} - \frac{1}{\overline{x}_i - x_{t,i}} \right),
     (\boldsymbol{\omega}_{\mathbf{u}})_j = -\sigma_j \left( \frac{1}{u_{t,j} - \underline{u}_j} - \frac{1}{\overline{u}_j - u_{t,j}} \right) \nonumber \\
    &(\boldsymbol{\omega}_{\mathbf{xx}})_{ii} = \mu_i \left( \frac{1}{(x_{t,i} - \underline{x}_i)^2} + \frac{1}{(\overline{x}_i - x_{t,i})^2} \right),\nonumber \\
     &(\boldsymbol{\omega}_{\mathbf{uu}})_{jj} = \sigma_j \left( \frac{1}{(u_{t,j} - \underline{u}_j)^2} + \frac{1}{(\overline{u}_j - u_{t,j})^2} \right),
\end{align}
where $i \in \mathcal{I}_{x}$ and $j \in \mathcal{I}_{u}$. The cross-term Hessian $\boldsymbol{\omega}_{\mathbf{ux}}$ is zero.\\ The \textbf{terminal barrier term} is defined as:\\ $\boldsymbol{\Omega}(\mathbf{x}_T) = - \sum_{i \in \mathcal{I}_x} \mu_{i}\left( \log(x_{T,i} - \underline{x}_{i}) + \log(\overline{x}_i - x_{T,i}) \right)$.
The non-zero elements of its derivatives with respect to the terminal state $\mathbf{x}_T$ are:\\ $ (\boldsymbol{\Omega}_{\mathbf{x}})_i = -\mu_i \left( \frac{1}{x_{T,i} - \underline{x}_i} - \frac{1}{\overline{x}_i - x_{T,i}} \right)$ and $(\boldsymbol{\Omega}_{\mathbf{xx}})_{ii} = \mu_i \left( \frac{1}{(x_{T,i} - \underline{x}_i)^2} + \frac{1}{(\overline{x}_i - x_{T,i})^2} \right)$.

\subsubsection{Backward Pass}
The backward pass computes the local optimal control policy by propagating the coefficients of the updated costate $\boldsymbol{\lambda}_{t}$ = $\mathbf{S}_t\delta \mathbf{x}_{t} + \mathbf{v}_{t}$. It computes the factors, $\mathbf{v}_t$ and $\mathbf{S}_t$, backward in time from $t=T$ to $t=0$.

\noindent\textbf{1. Initialization ($t=T$):}
The value function derivatives are initialized at the terminal time using the derivatives of the total terminal cost, evaluated at the nominal state $\bar{\mathbf{x}}_T$: $ \mathbf{v}_{T} = \Phi_{\mathbf{x}} + \boldsymbol{\Omega}_{\mathbf{x}} \quad \text{and} \quad \mathbf{S}_{T} = \Phi_{\mathbf{xx}} + \boldsymbol{\Omega}_{\mathbf{xx}}$.

\noindent\textbf{2. Backward Recursion ($t=T-1, \dots, 0$):}
For each time step, we compute the derivatives of the action-value function $Q_t$ for deviations around $(\bar{\mathbf{x}}_t, \bar{\mathbf{u}}_t)$:
\begin{align}
    Q_{\mathbf{x}} &= (C_{\mathbf{x}} + \boldsymbol{\omega}_{\mathbf{x}}) + \mathbf{F}_{\mathbf{x}}^{\top}\mathbf{v}_{t+1} \\
    Q_{\mathbf{u}} &= (C_{\mathbf{u}} + \boldsymbol{\omega}_{\mathbf{u}}) + \mathbf{F}_{\mathbf{u}}^{\top}\mathbf{v}_{t+1} \\
    Q_{\mathbf{xx}} &= (C_{\mathbf{xx}} + \boldsymbol{\omega}_{\mathbf{xx}}) + \mathbf{F}_{\mathbf{x}}^{\top}\mathbf{S}_{t+1}\mathbf{F}_{\mathbf{x}} \label{eq:Q_xx}\\
    Q_{\mathbf{uu}} &= (C_{\mathbf{uu}} + \boldsymbol{\omega}_{\mathbf{uu}}) + \mathbf{F}_{\mathbf{u}}^{\top}\mathbf{S}_{t+1}\mathbf{F}_{\mathbf{u}}
    \label{eq:Q_uu}\\
    Q_{\mathbf{ux}} &= C_{\mathbf{ux}} + \mathbf{F}_{\mathbf{u}}^{\top}\mathbf{S}_{t+1}\mathbf{F}_{\mathbf{x}}
\end{align}
All derivatives are evaluated at the nominal trajectory, $(\mathbf{\tilde{x}}_{t},\mathbf{\tilde{u}}_t)$. The feed-forward term $\mathbf{k}_t$ and feedback gain matrix $\mathbf{K}_t$ are then computed. A regularization term $\zeta \mathbf{I}$ is added to $Q_{\mathbf{uu}}$ to ensure it is positive-definite.
\begin{align}
    \mathbf{k}_{t} &= -(Q_{\mathbf{uu}} + \zeta\mathbf{I})^{-1}Q_{\mathbf{u}} \label{eq:k_update}\\
    \mathbf{K}_{t} &= -(Q_{\mathbf{uu}} + \zeta\mathbf{I})^{-1}Q_{\mathbf{ux}} \label{eq:K_update}
\end{align}
Finally, the values of $\mathbf{S}_t$ and $\mathbf{v}_t$ are updated:
\begin{align}
    \mathbf{v}_{t} &= Q_{\mathbf{x}} + \mathbf{K}_{t}^{\top}Q_{\mathbf{uu}}\mathbf{k}_{t} + \mathbf{K}_{t}^{\top}Q_{\mathbf{u}} + Q_{\mathbf{ux}}^{\top}\mathbf{k}_{t} \label{eq:v_update} \\
    \mathbf{S}_{t} &= Q_{\mathbf{xx}} + \mathbf{K}_{t}^{\top}Q_{\mathbf{uu}}\mathbf{K}_{t} + \mathbf{K}_{t}^{\top}Q_{\mathbf{ux}} + Q_{\mathbf{ux}}^{\top}\mathbf{K}_{t} \label{eq:S_update}
\end{align}

In the forward pass, a new trajectory is generated by applying the computed control policy. A line search is performed with a parameter $\alpha \in (0, 1]$ to ensure improvement in the cost function and maintain feasibility. Starting from $\mathbf{x}_0^{\text{new}} = \bar{\mathbf{x}}_0$:
\begin{align}
    \mathbf{u}_{t}^{\text{new}} &= \tilde{\mathbf{u}}_{t} + \alpha \mathbf{k}_{t} + \mathbf{K}_{t}(\mathbf{x}_{t}^{\text{new}} - \tilde{\mathbf{x}}_{t}) \quad & t \in [0, T-1] \label{eq:u_update}\\
    \mathbf{x}_{t+1}^{\text{new}} &= \mathbf{F}(\mathbf{x}_{t}^{\text{new}}, \mathbf{u}_{t}^{\text{new}}) \quad & t \in [0, T-1]
\end{align}
The new trajectory $(\mathbf{x}^{\text{new}}, \mathbf{u}^{\text{new}})$ becomes the nominal trajectory for the next iteration. This forward-backward procedure is repeated until convergence. The backtracking line search on $\alpha$ ensures that the augmented cost function $\mathcal{J}_{\boldsymbol{\mu},\boldsymbol{\sigma}}$ decreases at each iteration. A new trajectory is accepted if it provides a sufficient decrease, typically measured by comparing the actual reduction in cost to the reduction predicted by the local quadratic model. This ensures the algorithm makes progress towards a local minimum. A common acceptance criterion is when the ratio of actual to expected reduction is positive:
\begin{equation}
\label{eq:ilqr_reduction_ratio}
    \frac{\mathcal{J}(\text{new}) - \mathcal{J}(\text{old})}{\Delta J_{\text{expected}}} > c_1
\end{equation}
where $\mathcal{J}(\text{old})$ and $\mathcal{J}(\text{new})$ are the costs of the old and new trajectories, $c_1 \in (0, 1)$ is a small constant (e.g., $10^{-4}$), and $\Delta J_{\text{expected}} < 0$ is the cost change predicted by the quadratic approximation.
The predicted change in the local quadratic model is
\begin{equation}
    \Delta J_{\mathrm{expected}}
    =
    \sum_{t=0}^{T-1}
    \left(
        \alpha \mathbf{k}_t^{\top} Q_{\mathbf{u},t}
        +
        \frac{\alpha^2}{2}
        \mathbf{k}_t^{\top}
        Q_{\mathbf{uu},t}
        \mathbf{k}_t
    \right).
    \label{eq:expected_Cost_reduction1}
\end{equation}
The regularization parameter $\zeta$ is adjusted so that
$Q_{\mathbf{uu},t}+\zeta\mathbf{I}$ is positive definite, ensuring that
the computed feedforward direction is a descent direction. The subsequent
line search selects $\alpha$ to obtain an acceptable decrease in the
barrier-augmented objective.

\begin{algorithm}[H]
\caption{Box-iLQR}
\label{alg:barrier_update}
\begin{algorithmic}[1]
\Require Initial feasible trajectory $(\tilde{\mathbf{x}}, \tilde{\mathbf{u}})$
\Require Initial barrier parameters $\boldsymbol{\mu}_0, \boldsymbol{\sigma}_0$
\Require Initial reduction factors $\mathbf{r}_{\mu} \in (0,1)^{n_1}, \mathbf{r}_{\sigma} \in (0,1)^{m_1}$
\Require Reduction factor update rate $\beta_r > 1$ (e.g., 1.5)
\Require Barrier termination tolerance $\epsilon_{\text{barrier}}$

\State Initialize $\boldsymbol{\mu} \leftarrow \boldsymbol{\mu}_0$, $\boldsymbol{\sigma} \leftarrow \boldsymbol{\sigma}_0$
\State Initialize $\boldsymbol{\mu}_{\text{prev}} \leftarrow \boldsymbol{\mu}$, $\boldsymbol{\sigma}_{\text{prev}} \leftarrow \boldsymbol{\sigma}$

\While{$\|\boldsymbol{\mu}\| > \epsilon_{\text{barrier}}$ or $\|\boldsymbol{\sigma}\| > \epsilon_{\text{barrier}}$}
    \State \Comment{Solve the inner-loop problem with a warm start.}
    \State $(\bar{\mathbf{x}}, \bar{\mathbf{u}}), \text{success}, \text{failed\_idx} \leftarrow \Call{iLQR\_Solve}{\bar{\mathbf{x}}, \bar{\mathbf{u}}, \boldsymbol{\mu}, \boldsymbol{\sigma}}$
    
    \If{\textbf{not} success}
        \State \Comment{Reduction was too aggressive; revert and adapt.}
        \State $\boldsymbol{\mu} \leftarrow \boldsymbol{\mu}_{\text{prev}}$, $\boldsymbol{\sigma} \leftarrow \boldsymbol{\sigma}_{\text{prev}}$
        \If{constraint $i$ in $\boldsymbol{\mu}$ failed (from $\text{failed\_idx}$)}
            \State $(\mathbf{r}_{\mu})_i \leftarrow \min(1, (\mathbf{r}_{\mu})_i \cdot \beta_r)$ \Comment{Make reduction less aggressive.}
        \ElsIf{constraint $j$ in $\boldsymbol{\sigma}$ failed}
            \State $(\mathbf{r}_{\sigma})_j \leftarrow \min(1, (\mathbf{r}_{\sigma})_j \cdot \beta_r)$
        \EndIf
        \State \textbf{continue} \Comment{Retry solving with old parameters and new rates.}
    \Else
        \State \Comment{Successful solve; prepare for next reduction.}
        \State $\boldsymbol{\mu}_{\text{prev}} \leftarrow \boldsymbol{\mu}$, $\boldsymbol{\sigma}_{\text{prev}} \leftarrow \boldsymbol{\sigma}$
        \State $\boldsymbol{\mu} \leftarrow \mathbf{r}_{\mu} \odot \boldsymbol{\mu}$ \Comment{Element-wise multiplication.}
        \State $\boldsymbol{\sigma} \leftarrow \mathbf{r}_{\sigma} \odot \boldsymbol{\sigma}$
    \EndIf
\EndWhile

\State \Return Optimal trajectory $(\tilde{\mathbf{x}}, \tilde{\mathbf{u}})$
\end{algorithmic}
\end{algorithm}

The Box-iLQR method presented in the previous section solves the optimization problem for a \textit{fixed} set of barrier parameters $\boldsymbol{\mu}$ and $\boldsymbol{\sigma}$. To recover the solution to the original constrained problem, these parameters must be driven to zero. This is managed by an outer loop that iteratively reduces the barrier parameters and uses the iLQR solver to find the solution for each new subproblem. This sequence of solutions ideally follows the \textit{central path} towards the constrained optimum. The complete procedure, which we refer to as the Box-iLQR algorithm, is detailed in Algorithm~\ref{alg:barrier_update}. The process begins with conservatively large initial barrier parameters, $\boldsymbol{\mu}_0$ and $\boldsymbol{\sigma}_0$, and a feasible initial trajectory. The algorithm then enters an inner loop where it first solves the iLQR subproblem for the current parameters, using the solution from the previous iteration as a warm start to accelerate convergence. 

A crucial aspect of this process is the parameter update rule. After a successful inner-loop solve, each component of $\boldsymbol{\mu}$ and $\boldsymbol{\sigma}$ is decreased by a multiplicative factor. A naive or overly aggressive reduction can push the warm-start trajectory into an infeasible region where the logarithmic barriers are undefined, causing the inner-loop solver to fail. To handle this, our update schedule is adaptive. If the iLQR solver fails to converge, it indicates that the barrier parameters were reduced too severely. In response, the algorithm reverts the parameter update and applies a more conservative reduction factor for the corresponding constraint in the next attempt. This process continues until the barrier parameters are below a specified tolerance, ensuring their contribution to the total cost is negligible and the solution closely approximates the true constrained optimum.


\subsection{Feedback Stabilization}
\label{subsec:feedback_stabilization}

In addition to generating a nominal trajectory, the iLQR backward pass provides
a time-varying feedback gain matrix for every barrier subproblem. Let
\[
\left\{
\bar{\mathbf{x}}_k^{\mathrm{fb}},
\bar{\mathbf{u}}_k^{\mathrm{fb}},
\mathbf{K}_k^{\mathrm{fb}}
\right\}_{k=0}^{N-1}
\]
denote the converged trajectory and feedback gains associated with a selected
finite barrier value $\sigma_{\mathrm{fb}}>0$. The corresponding closed-loop
controller is
\begin{equation}
\mathbf{u}_k^{\mathrm{cl}}
=
\Pi_{\mathcal U}
\left[
\bar{\mathbf{u}}_k^{\mathrm{fb}}
+
\mathbf{K}_k^{\mathrm{fb}}
\left(
\mathbf{x}_k-\bar{\mathbf{x}}_k^{\mathrm{fb}}
\right)
\right],
\label{eq:feedback_stabilization}
\end{equation}
where $\Pi_{\mathcal U}$ denotes projection onto the admissible control set.

Barrier continuation may subsequently proceed toward zero to obtain a nominal
trajectory that more closely approaches the constrained optimum. The feedback
policy used for closed-loop execution, however, need not be taken from this
smallest-barrier solution. Instead, a finite barrier can be retained to
preserve constraint-dependent curvature and available control authority.

Retaining a finite barrier preserves the constraint-dependent curvature in the local quadratic model and therefore produces a constraint-aware feedback policy. As a control approaches its saturation boundary, the barrier curvature increases, reducing feedback corrections through that channel. Similarly, the effect of an approaching state constraint is propagated backward through the value-function recursion and influences the preceding feedback gains. Thus, the barrier parameter serves two distinct purposes: it is reduced toward zero to recover the constrained optimal trajectory, while a finite value is retained to obtain a practical feedback controller for disturbance rejection and trajectory stabilization. A convergence analysis of Box-iLQR and a detailed characterization of the
finite-barrier feedback policy are provided in our previous work
\cite{abhijeet2026safeoptimalcontrolusing}.

\section{Numerical Results}
\label{sec:results}

Box-iLQR is evaluated on three nonlinear spacecraft optimal-control problems:
a fully actuated attitude maneuver, a planar minimum-fuel orbit transfer, and
a fixed-time $L_2$-to-$L_2$ Halo-orbit transfer in the Earth--Moon CR3BP.
Barrier continuation is used to approach the constrained nominal solution,
while feedback gains from selected finite-barrier subproblems are used for
closed-loop execution according to Eq.~\eqref{eq:feedback_stabilization}.

\subsection{Spacecraft Attitude Control}
\label{subsec:attitude_results}

The attitude state is
\begin{equation}
\mathbf{x}
=
\begin{bmatrix}
\phi & \theta & \psi &
\omega_1 & \omega_2 & \omega_3
\end{bmatrix}^{T},
\end{equation}
where $(\phi,\theta,\psi)$ are $3$--$2$--$1$ Euler angles and
$\boldsymbol{\omega}=[\omega_1,\omega_2,\omega_3]^T$ is the body-frame angular velocity. The dynamics are
\begin{align}
\begin{bmatrix}
\dot{\phi}\\
\dot{\theta}\\
\dot{\psi}
\end{bmatrix}
&=
\begin{bmatrix}
1 & \sin\phi\tan\theta & \cos\phi\tan\theta\\
0 & \cos\phi & -\sin\phi\\
0 & \dfrac{\sin\phi}{\cos\theta} &
\dfrac{\cos\phi}{\cos\theta}
\end{bmatrix}
\boldsymbol{\omega},\\
\dot{\boldsymbol{\omega}}
&=
\mathbf{I}^{-1}
\left(
\mathbf{u}
-
\boldsymbol{\omega}\times\mathbf{I}\boldsymbol{\omega}
\right),
\label{eq:attitude_dynamics_results}
\end{align}
with
\[
\mathbf{I}=\operatorname{diag}(140,120,130),
\qquad
-10^{-3}\leq u_i\leq10^{-3},
\quad i=1,2,3.
\]

The desired terminal condition is
$
\boldsymbol{\eta}_f=
\begin{bmatrix}
\pi & -\pi/4 & -\pi/2
\end{bmatrix}^{T},
~~
\boldsymbol{\omega}_f=\mathbf{0},
$
where $\boldsymbol{\eta}=[\phi,\theta,\psi]^T$. The maneuver is propagated using zero-order-hold controls and fourth-order Runge--Kutta integration with $\Delta t=1~\mathrm{s}$ and $t_f=3000~\mathrm{s}$, giving $N=3000$ control intervals.

The objective is
\begin{equation}
J=
\frac{1}{2}\mathbf{e}_N^T\mathbf{Q}_f\mathbf{e}_N
+
\sum_{k=0}^{N-1}
\left[
\frac{1}{2}\mathbf{e}_k^T\mathbf{Q}\mathbf{e}_k
+
\frac{1}{2}\mathbf{u}_k^T\mathbf{R}\mathbf{u}_k
\right]\Delta t,
\label{eq:attitude_discrete_cost}
\end{equation}
with $\mathbf{Q}=\mathbf{R}=5\mathbf{I}$ and
$\mathbf{Q}_f=50000\mathbf{I}$. The initial condition is $\mathbf{x}_{0} = \boldsymbol{0}_{6}$. The feedback gains used in the robustness study are retained at
$\sigma_{\mathrm{fb}} \approx 10^{-5}$.

\begin{figure*}[t]
\centering
\includegraphics[width=0.32\textwidth]{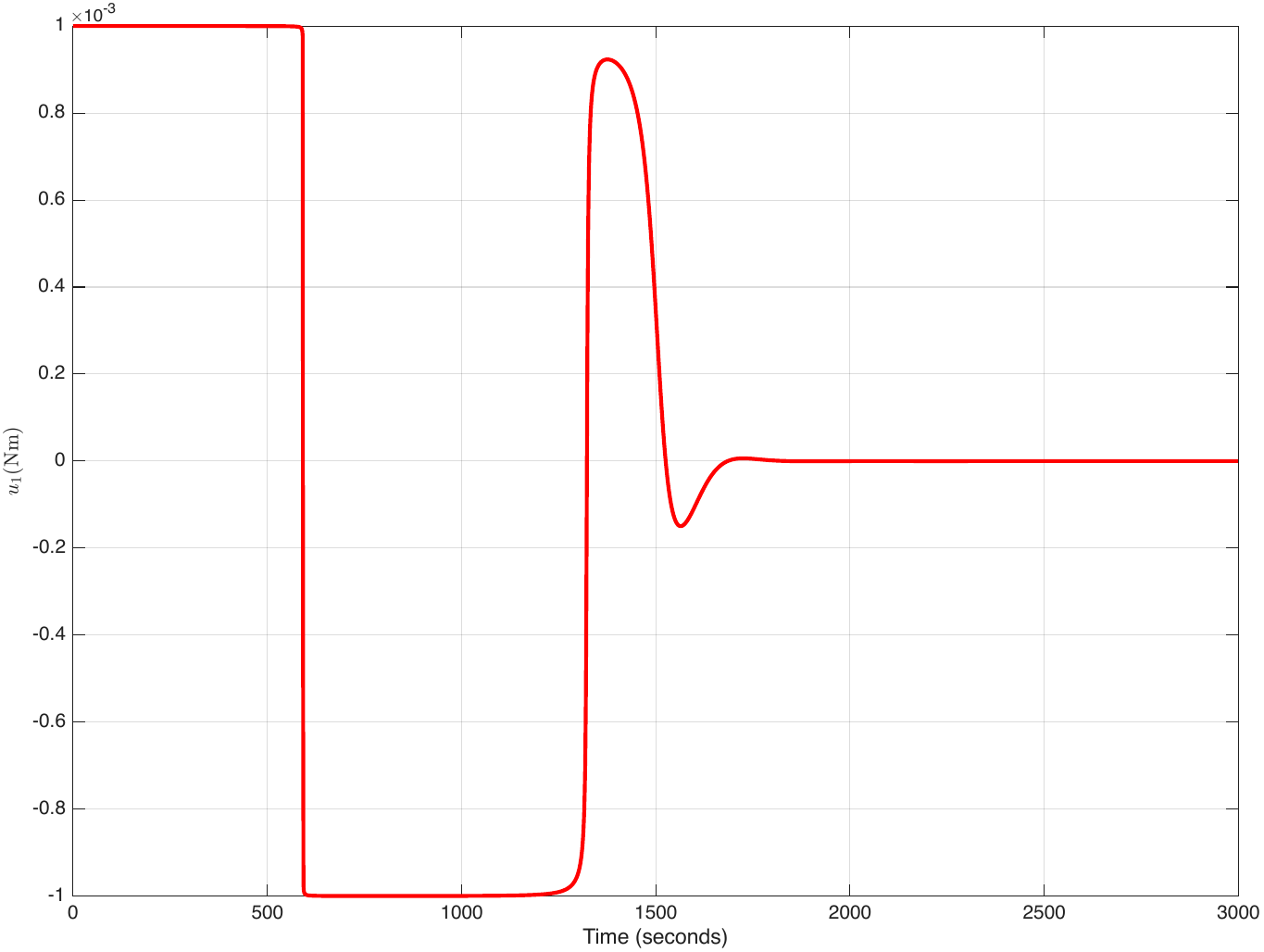}
\includegraphics[width=0.32\textwidth]{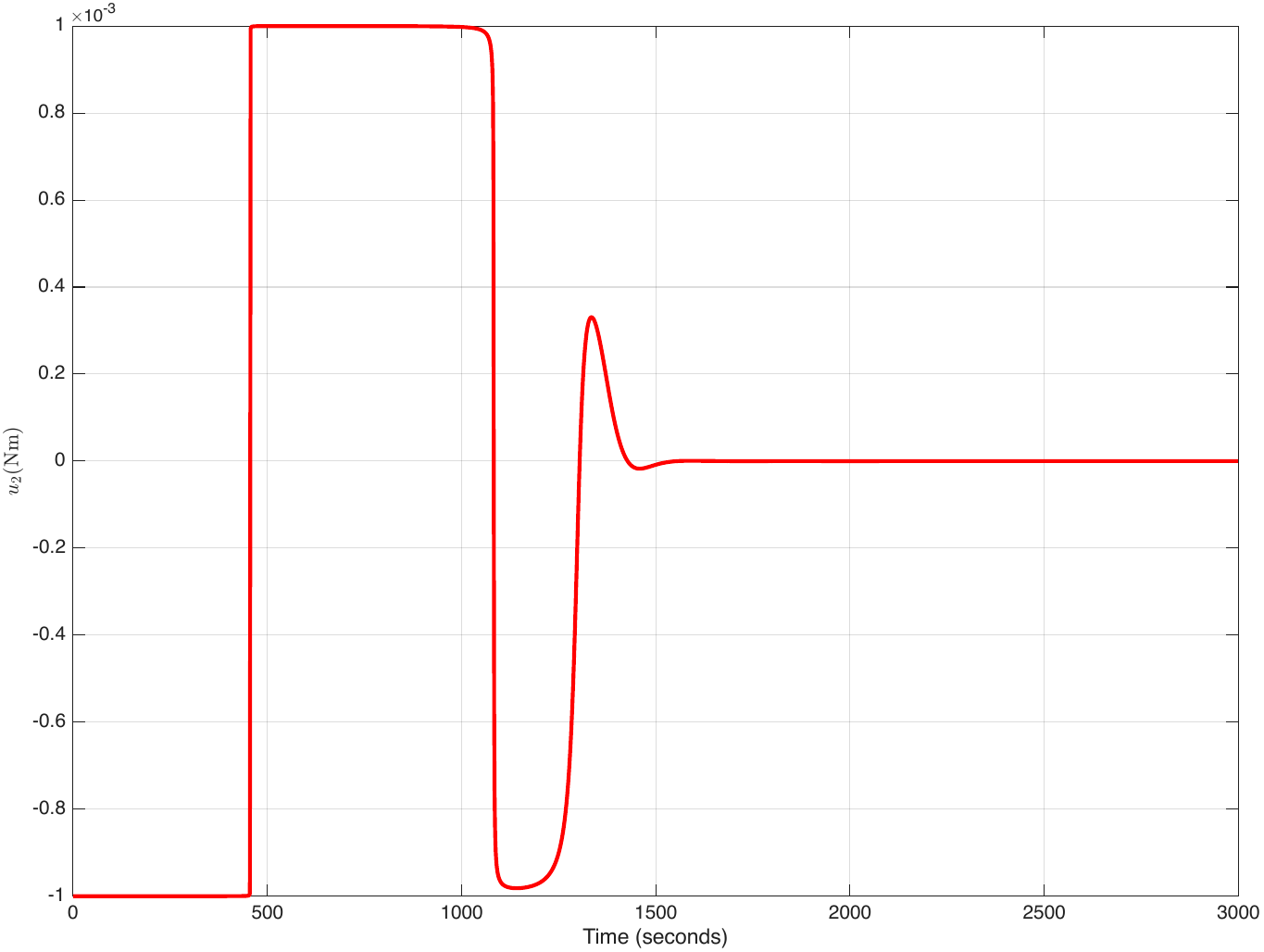}
\includegraphics[width=0.32\textwidth]{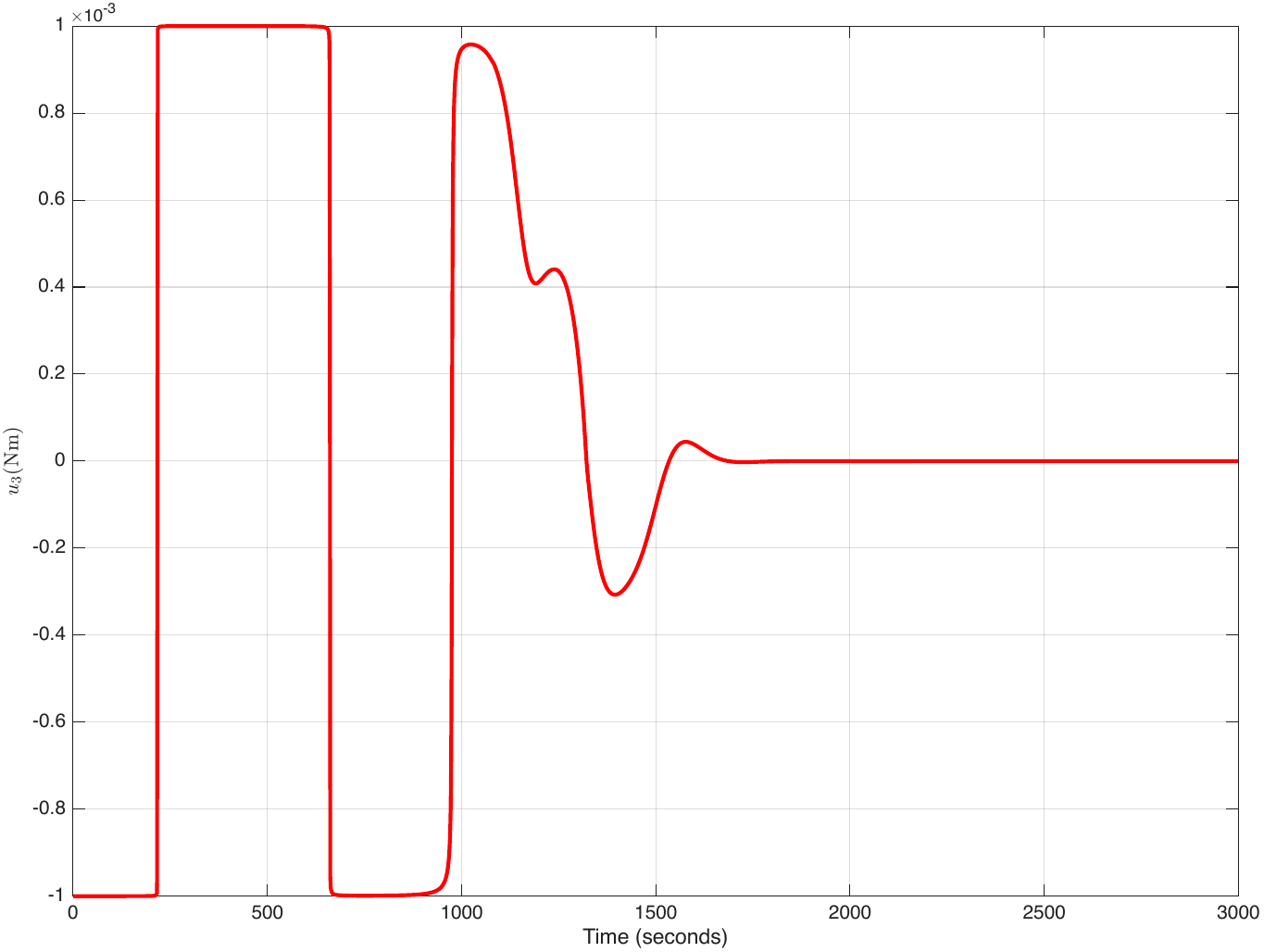}
\caption{Optimal attitude-control torque histories obtained using Box-iLQR. The torque histories satisfy $|u_i|\leq10^{-3}$ and contain extended boundary arcs.}
\label{fig:attitude_controls}
\end{figure*}

Figure~\ref{fig:attitude_controls} shows the optimal torque histories. The first torque remains at its upper bound for approximately $600~\mathrm{s}$, at its lower bound until approximately $1300~\mathrm{s}$, and settles near zero by $1800~\mathrm{s}$. The second remains at its lower bound for approximately $450~\mathrm{s}$, switches to its upper bound until approximately $1000~\mathrm{s}$, and approaches zero by $1500~\mathrm{s}$. The third remains initially at its lower bound for approximately $200~\mathrm{s}$, switches to its upper bound until approximately $700~\mathrm{s}$, briefly returns to its lower bound, and approaches zero by approximately $1700~\mathrm{s}$. Thus, the optimal maneuver contains extended control-boundary arcs.

\begin{figure*}[t]
\centering
\begin{minipage}[t]{0.49\textwidth}
\centering
\includegraphics[width=\linewidth]{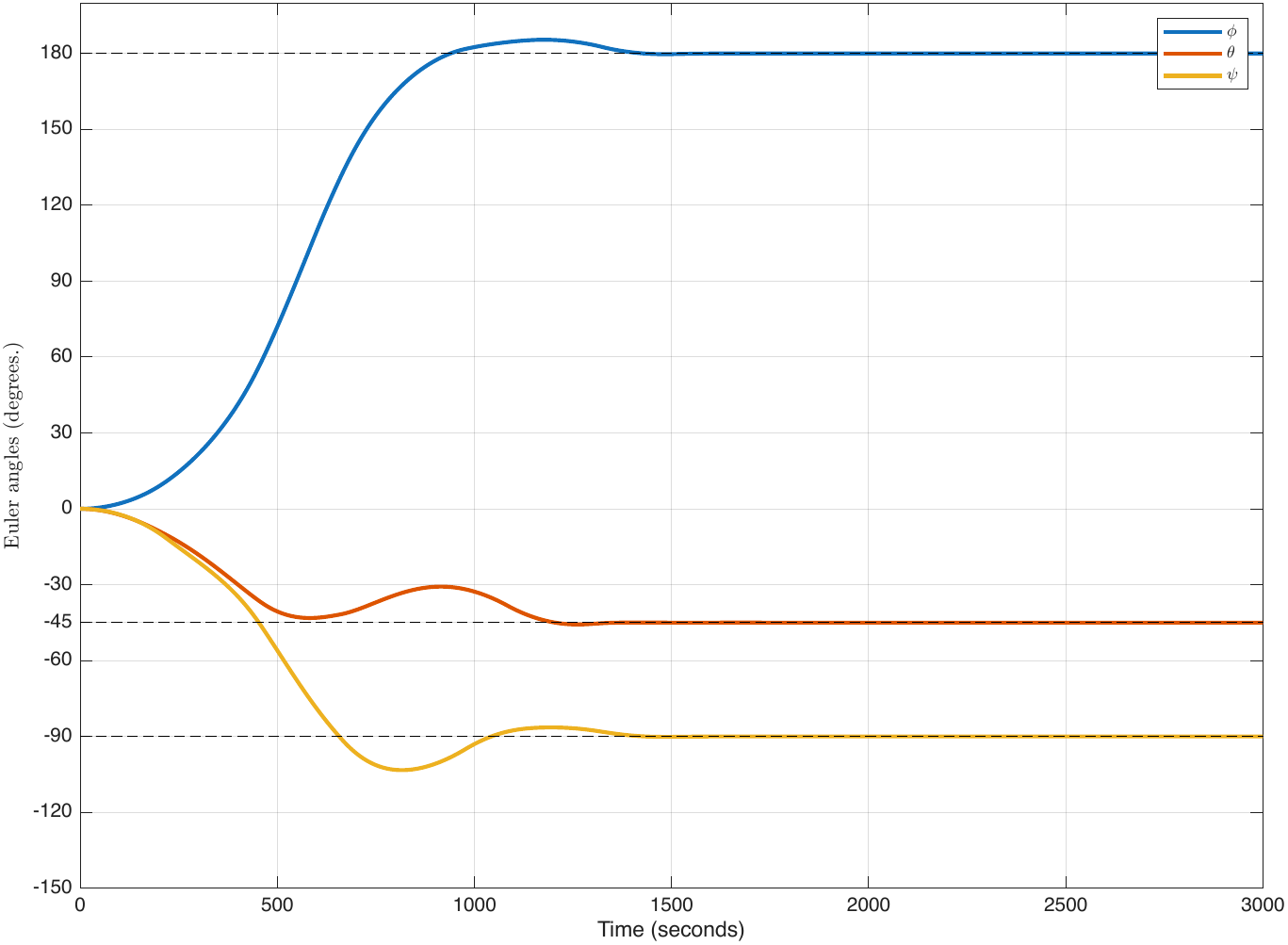}
(a)
\end{minipage}
\hfill
\begin{minipage}[t]{0.49\textwidth}
\centering
\includegraphics[width=\linewidth]{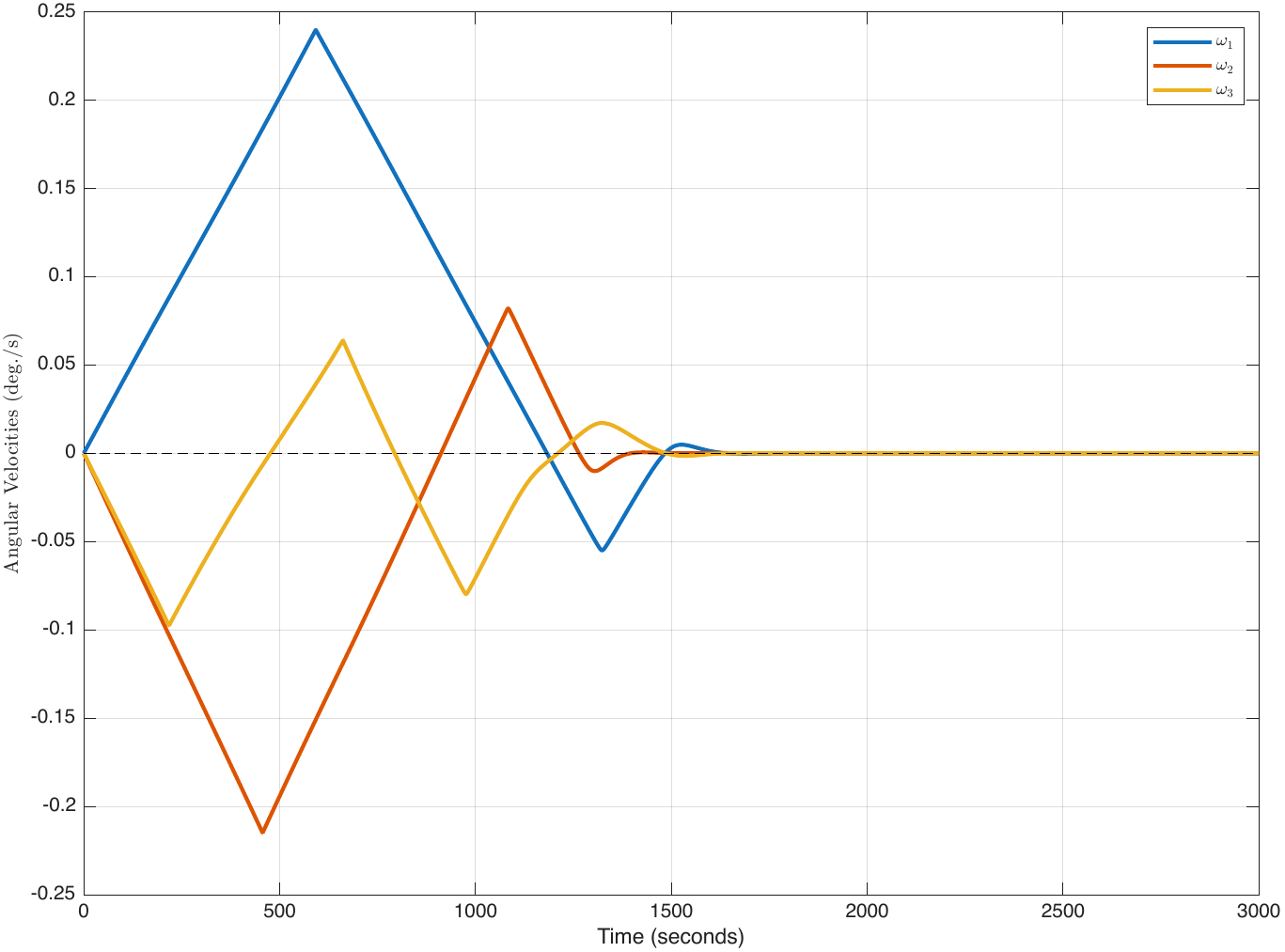}
(b)
\end{minipage}
\caption{Nominal attitude trajectory: (a) Euler angles and (b) body-frame angular velocities.}
\label{fig:attitude_states_rates}
\end{figure*}

As shown in Fig.~\ref{fig:attitude_states_rates}, the Euler angles converge to $\pi$, $-\pi/4$, and $-\pi/2$ before approximately $1500~\mathrm{s}$. The angular velocities exhibit approximately linear growth and decay during the saturated-control intervals and subsequently converge to zero, with the dominant transient behavior disappearing by approximately $1600~\mathrm{s}$.

Robustness is evaluated using simultaneous disturbance torques and measurement noise. The disturbance is
\begin{equation}
\boldsymbol{\tau}_d(t)=
\boldsymbol{\xi}
+
\mathbf{a}\sin(\boldsymbol{\nu}t)
+
\mathbf{b}\cos(\boldsymbol{\nu}t),
\label{eq:attitude_disturbance}
\end{equation}
with independently sampled components

$$
\xi_i\sim\mathcal{N}(0,(10^{-5})^2),\qquad
a_i,b_i\sim\mathcal{N}(0,(10^{-6})^2),\qquad
\nu_i\sim\mathcal{N}(0,(10^{-4})^2).
$$

The perturbed rotational dynamics are therefore

$$
\dot{\boldsymbol{\omega}}
=
\mathbf{I}^{-1}
\left[
\mathbf{u}
+\boldsymbol{\tau}_d(t)
-\boldsymbol{\omega}\times\mathbf{I}\boldsymbol{\omega}
\right].
$$

Measurement noise is modeled as
\begin{align*}
\widehat{\boldsymbol{\eta}}_k
=
\boldsymbol{\eta}_k+\mathbf{n}_{\eta,k},\quad
\mathbf{n}_{\eta,k}
\sim
\mathcal{N}
\left(
\mathbf{0},(10^{-4})^2\mathbf{I}_3
\right),\\
\widehat{\boldsymbol{\omega}}_k
=
\boldsymbol{\omega}_k+\mathbf{n}_{\omega,k},\quad
\mathbf{n}_{\omega,k}
\sim
\mathcal{N}
\left(
\mathbf{0},(10^{-6})^2\mathbf{I}_3
\right).
\end{align*}

\begin{figure*}[t]
    \centering
    \includegraphics[width=0.32\textwidth]{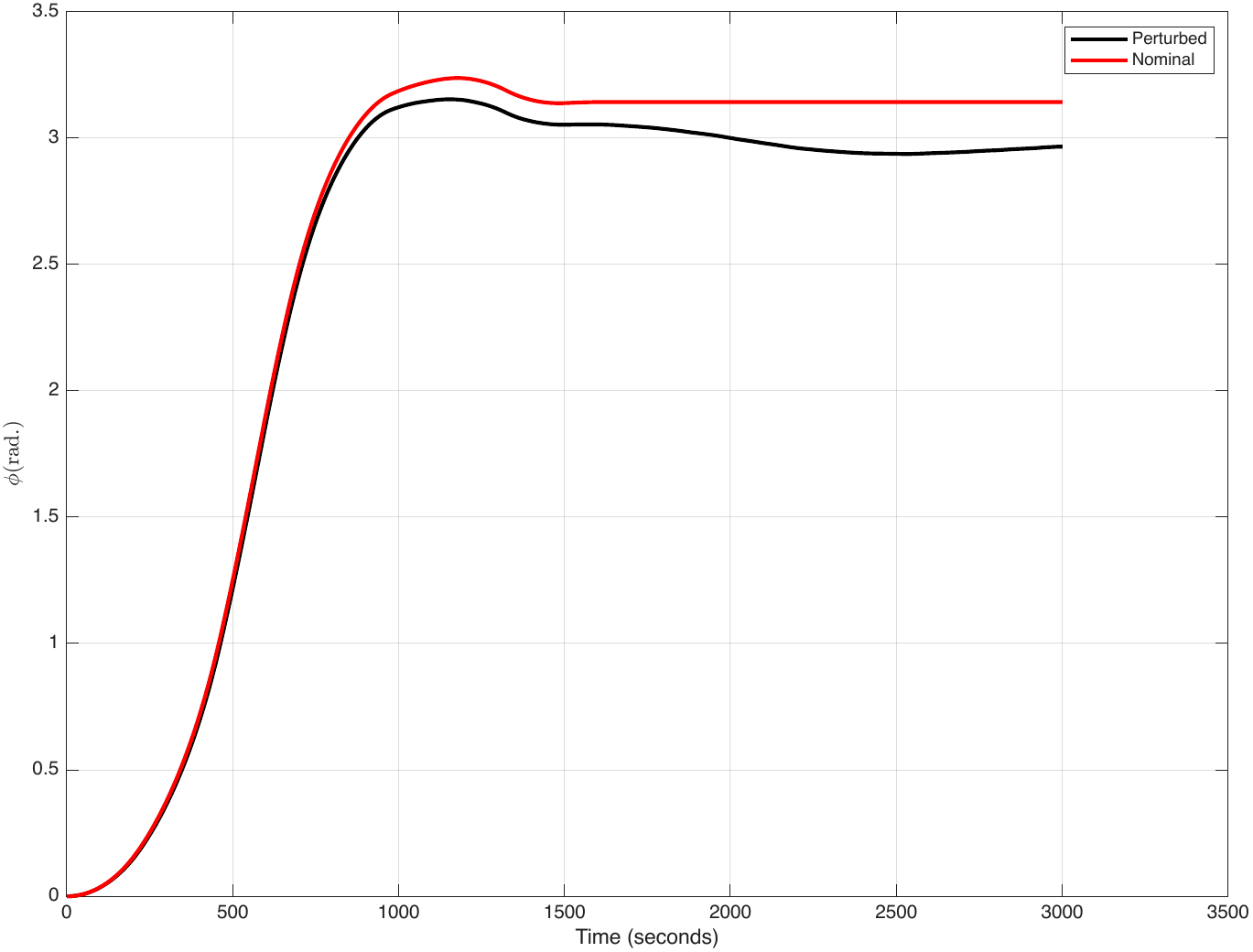}
    \includegraphics[width=0.32\textwidth]{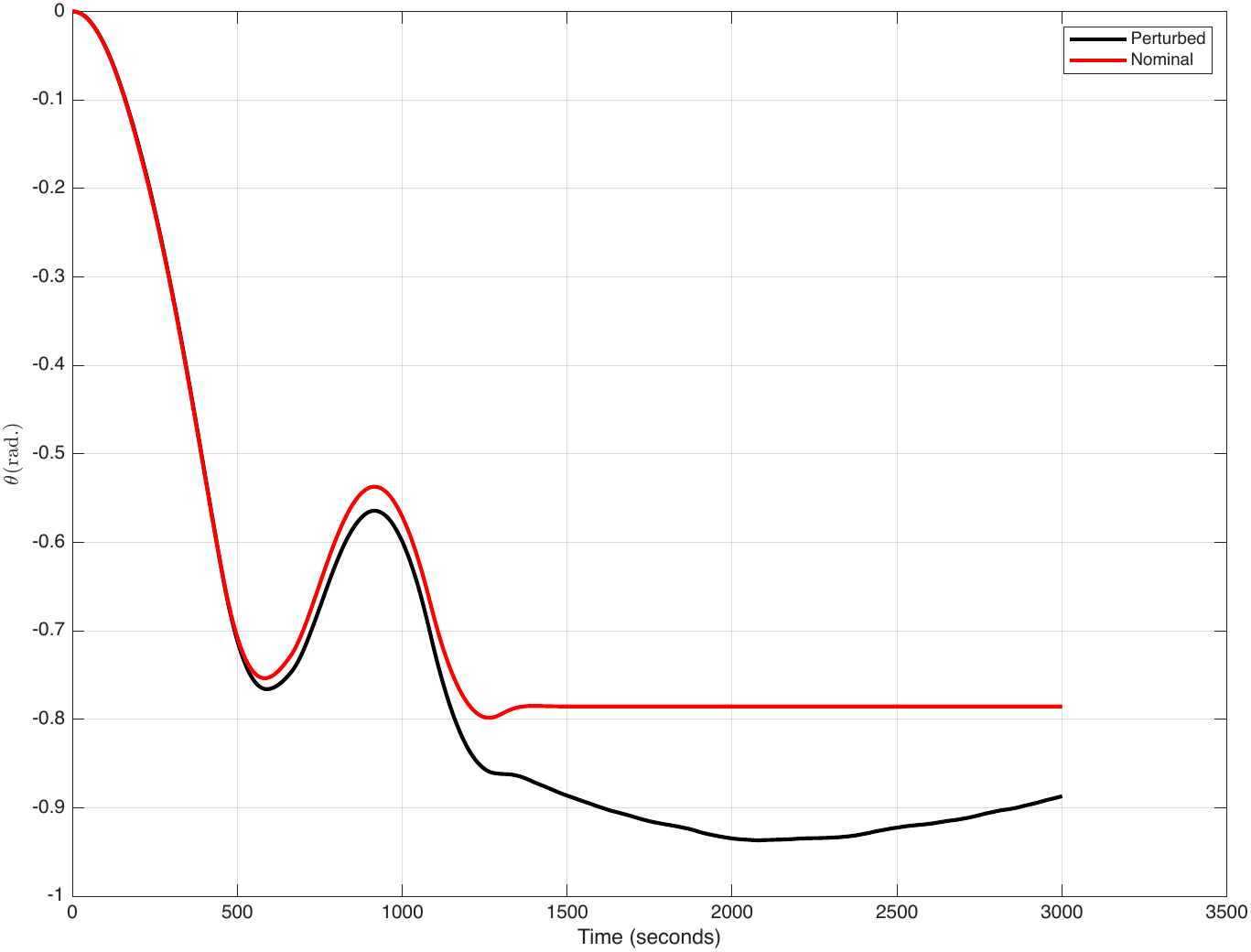}
    \includegraphics[width=0.32\textwidth]{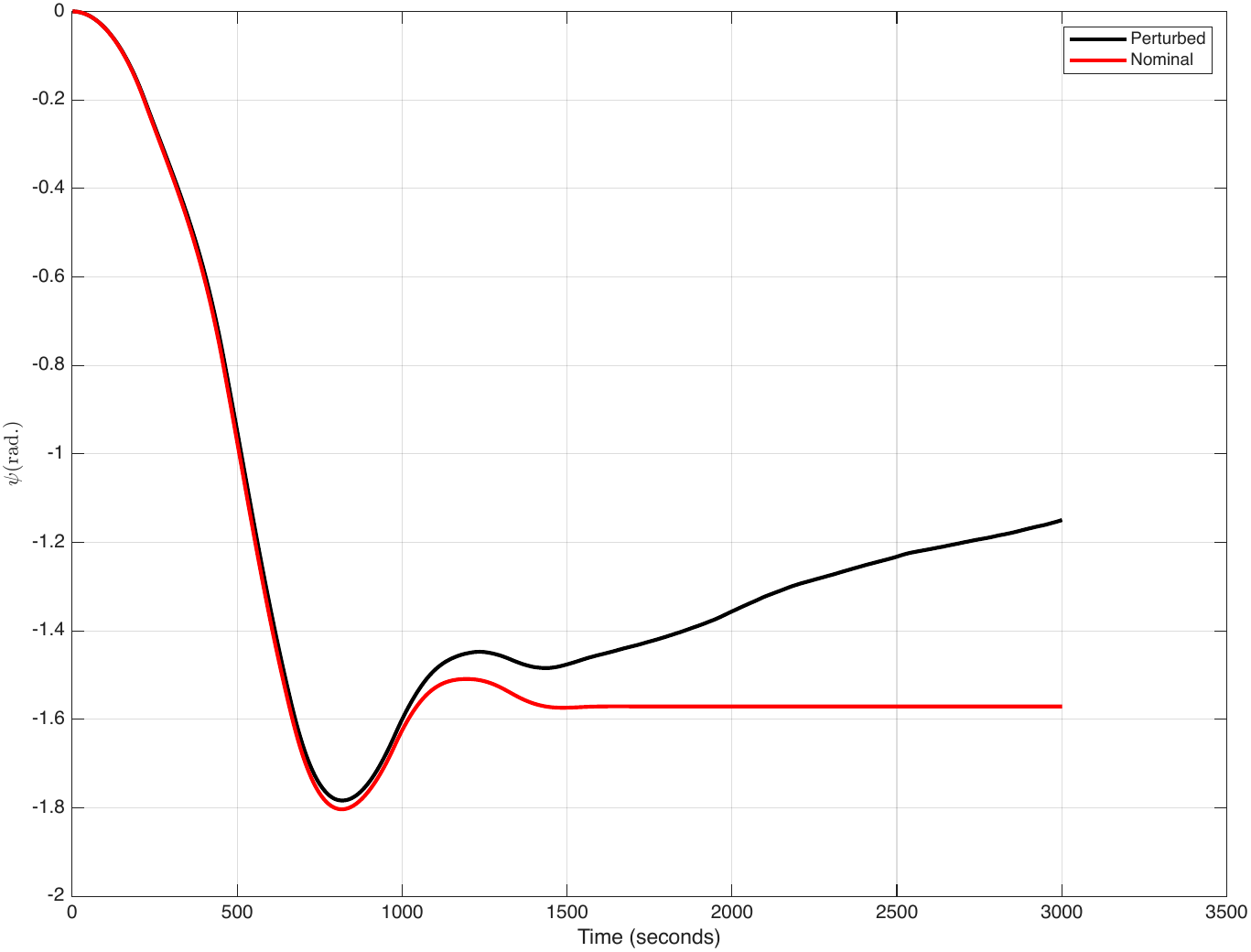}\\
    \includegraphics[width=0.32\textwidth]{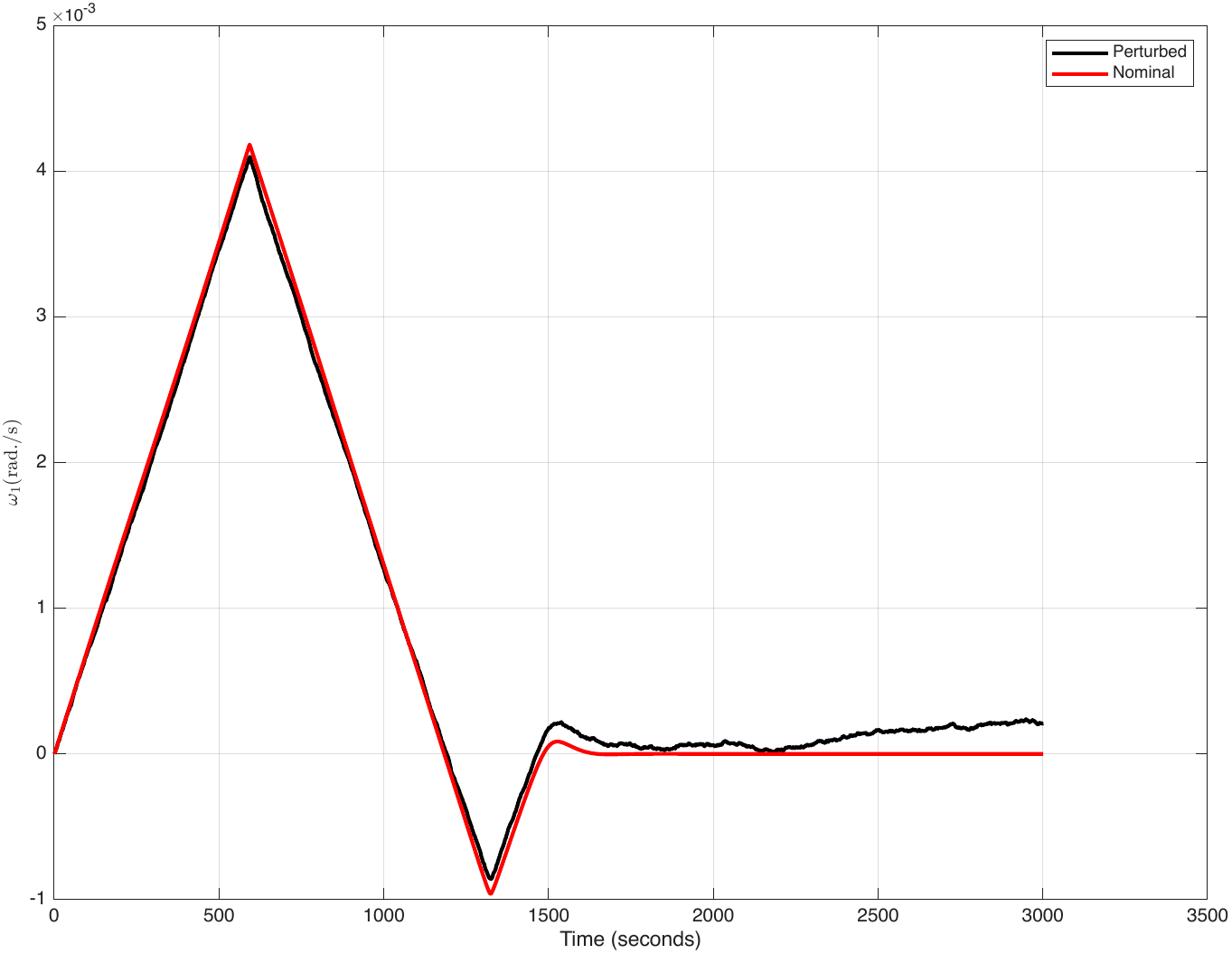}
    \includegraphics[width=0.32\textwidth]{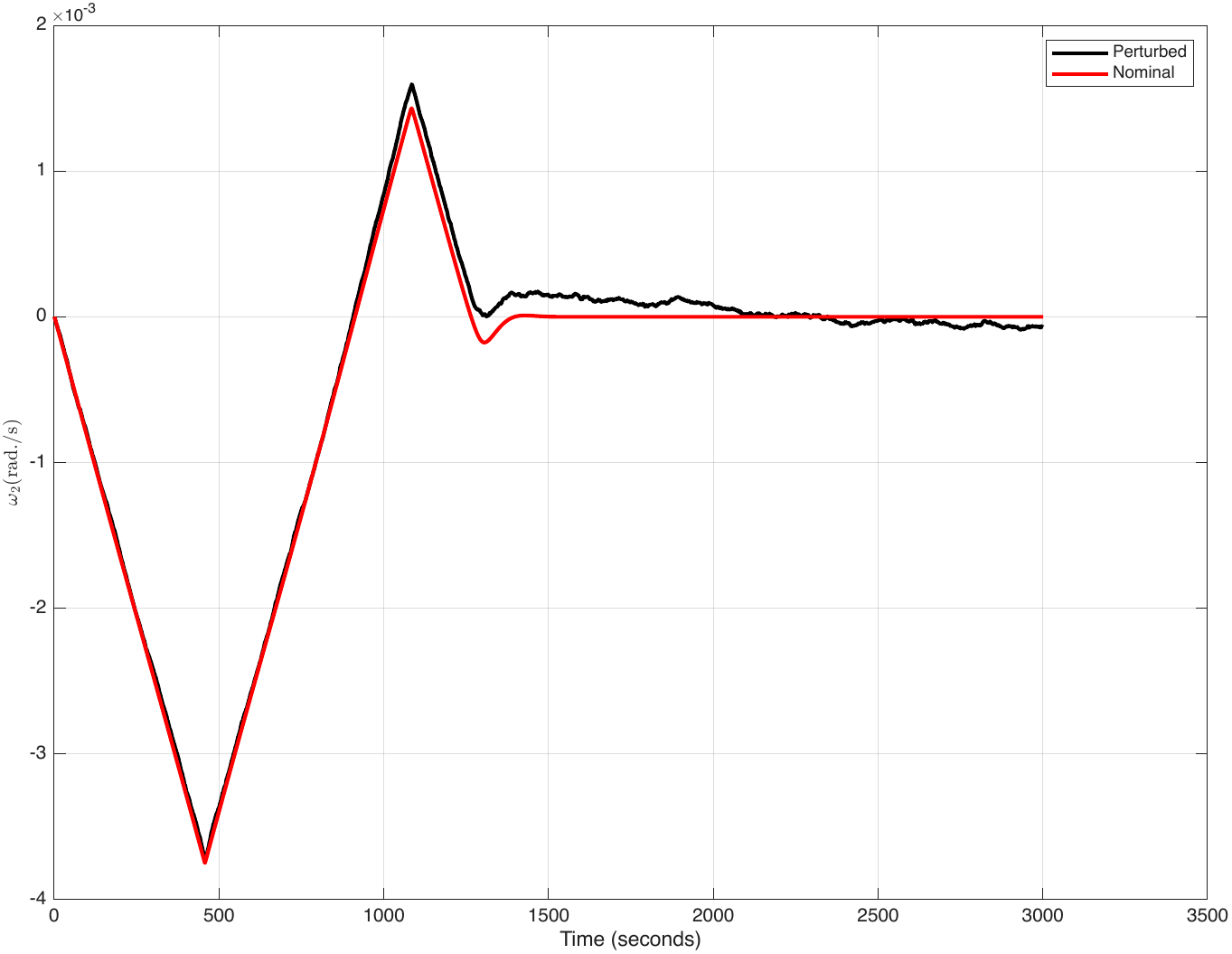}
    \includegraphics[width=0.32\textwidth]{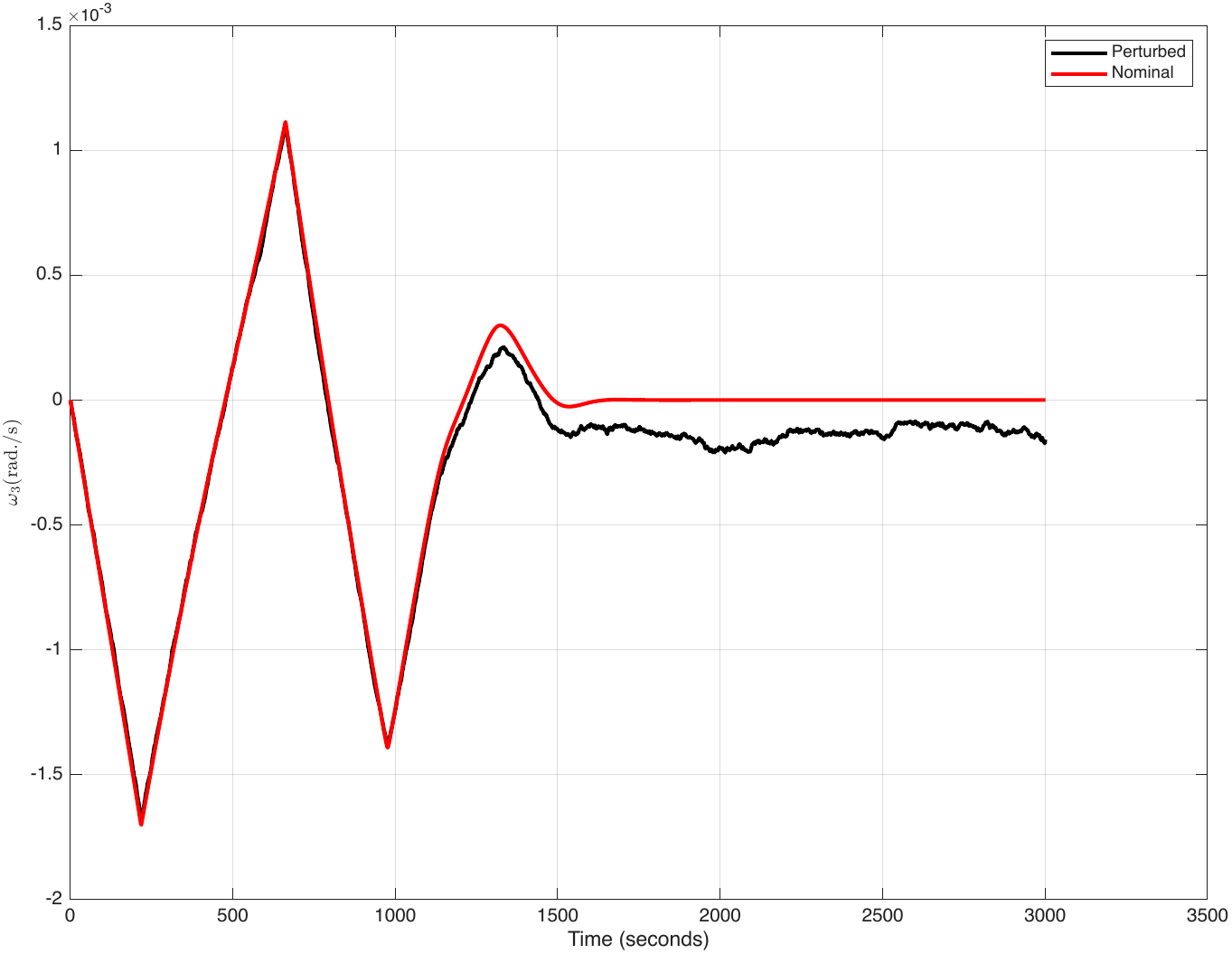}
    \caption{Attitude response under disturbances and measurement noise without feedback. The nominal trajectory is shown in red and the disturbed open-loop trajectory in black.}
    \label{fig:attitude_open_loop_perturbed}
\end{figure*}

\begin{figure*}[t]
    \centering
    \includegraphics[width=0.32\textwidth]{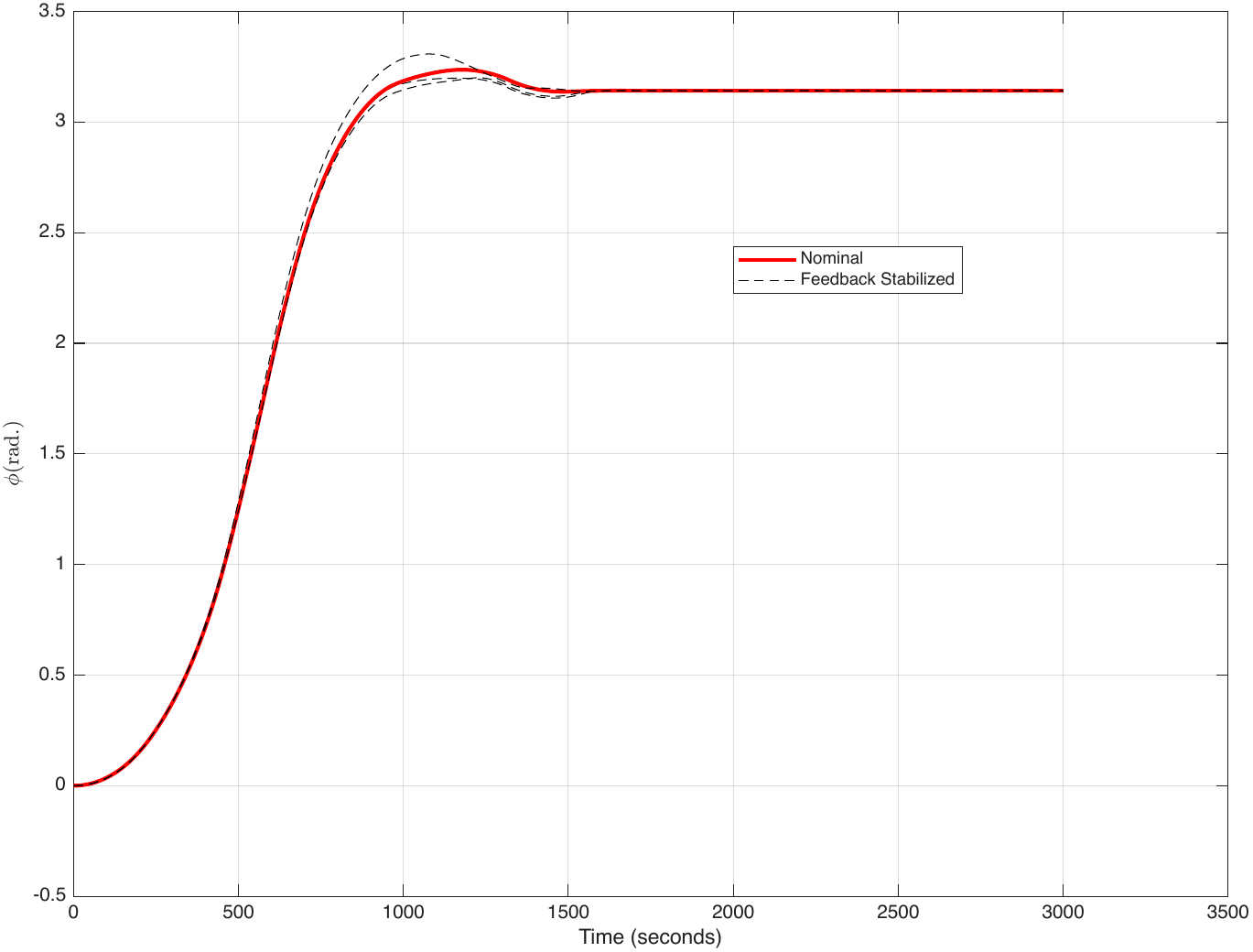}
    \includegraphics[width=0.32\textwidth]{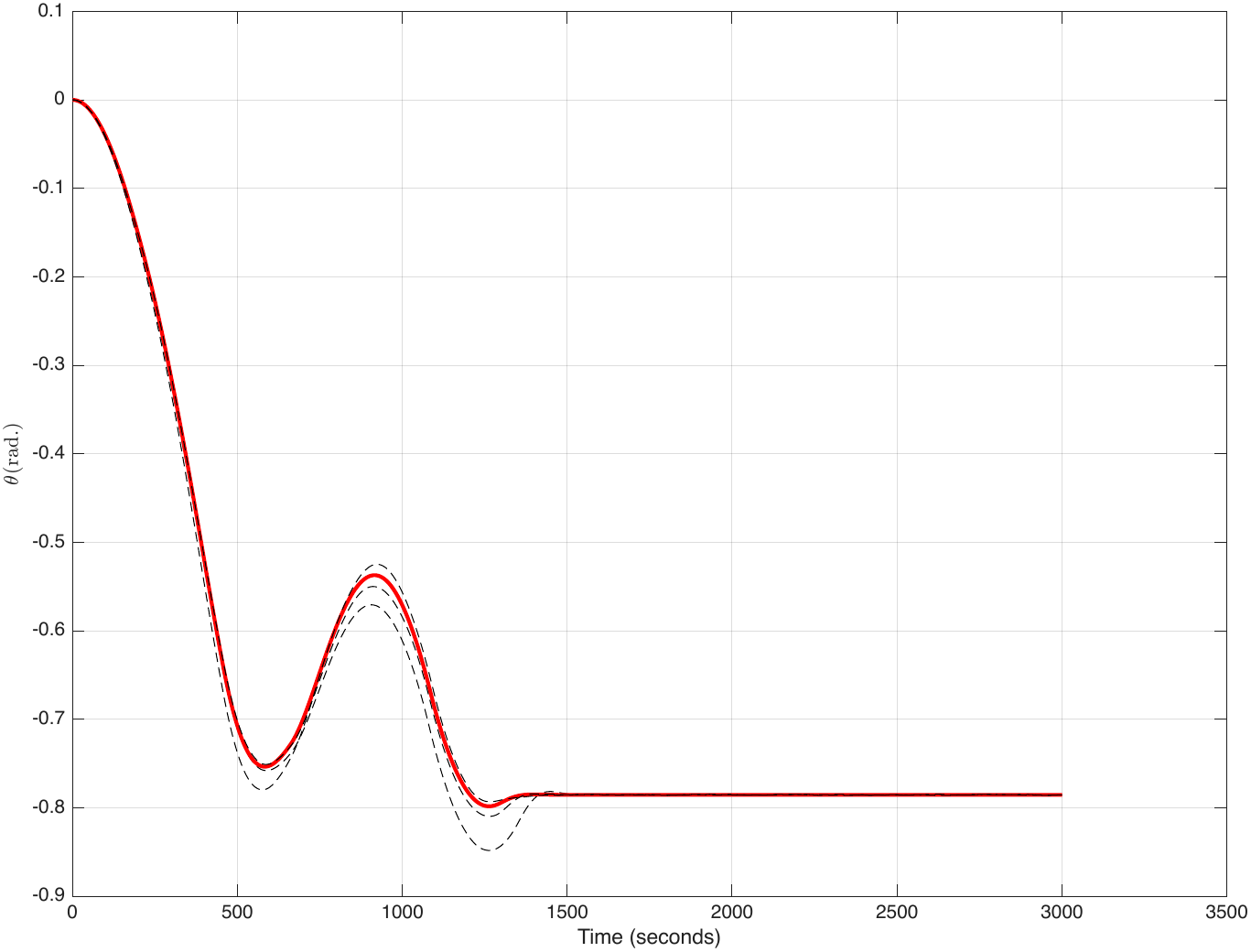}
    \includegraphics[width=0.32\textwidth]{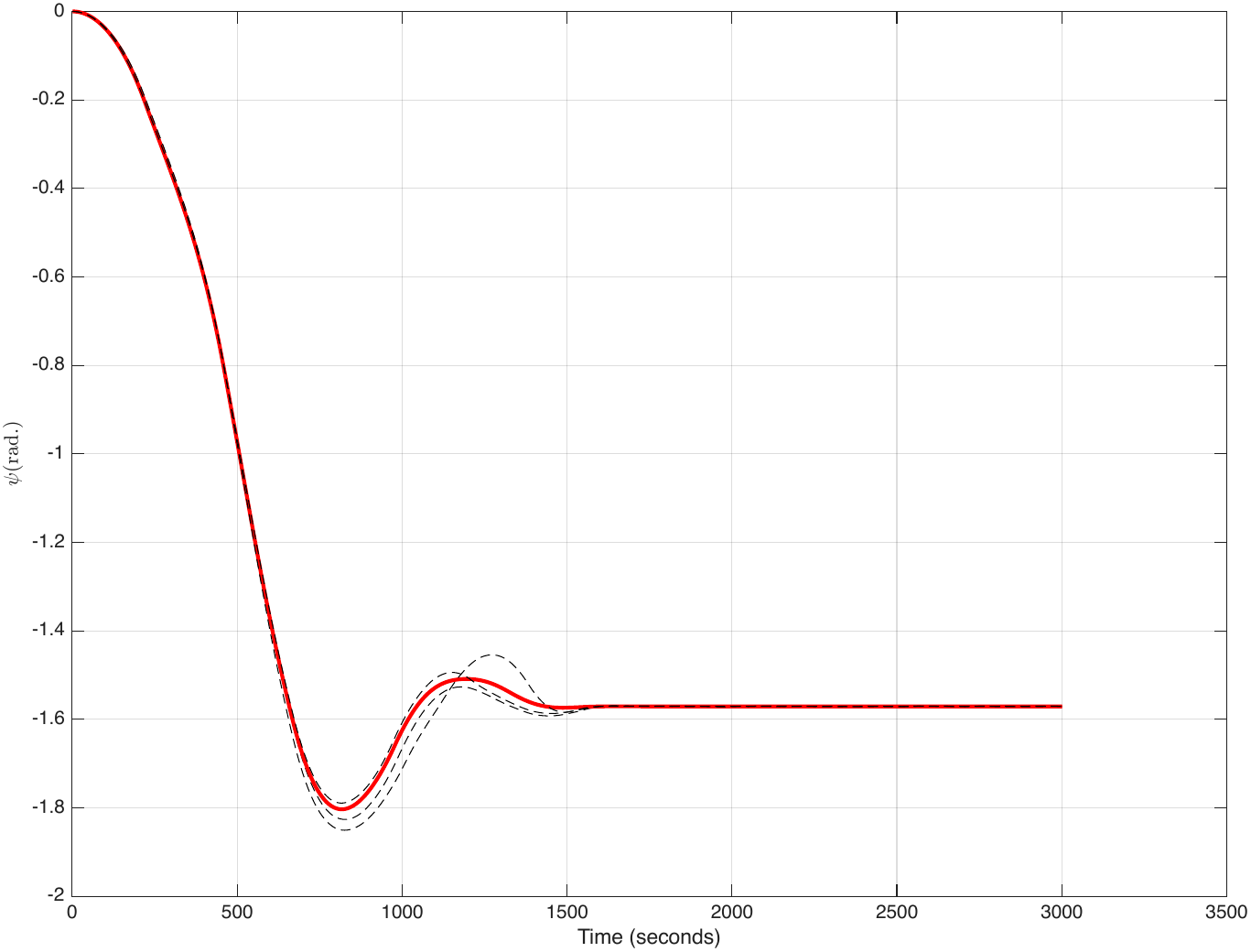}\\
    \includegraphics[width=0.32\textwidth]{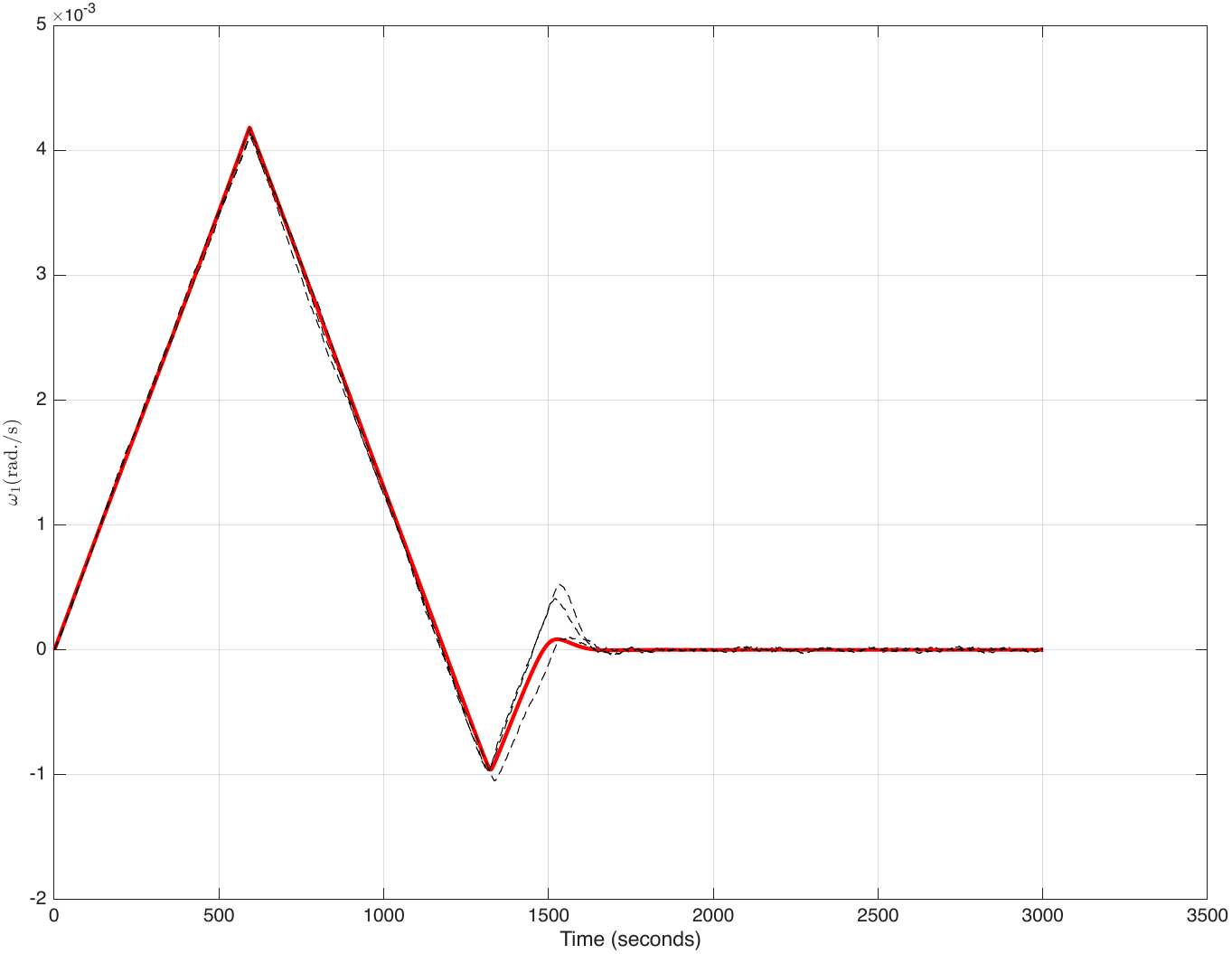}
    \includegraphics[width=0.32\textwidth]{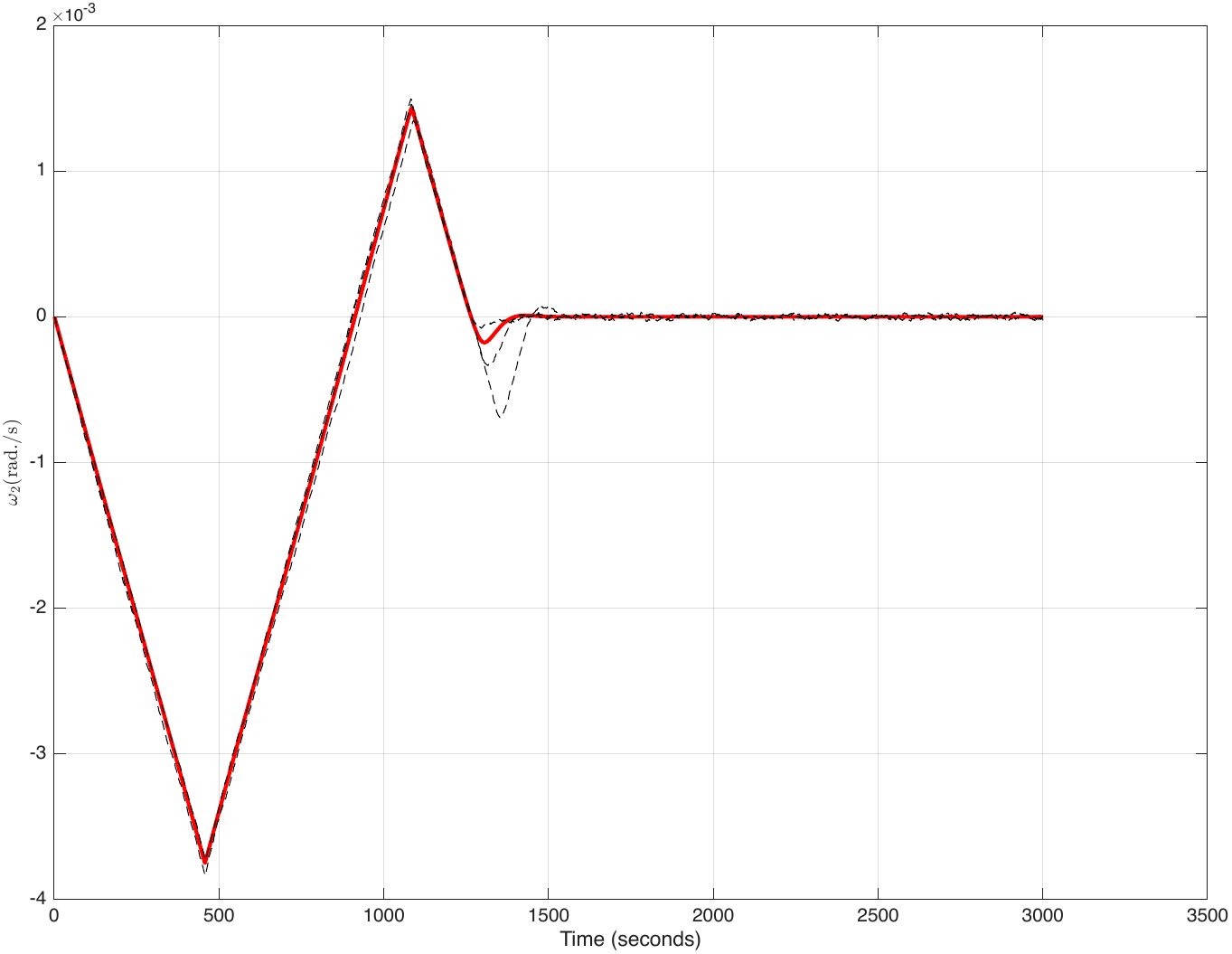}
    \includegraphics[width=0.32\textwidth]{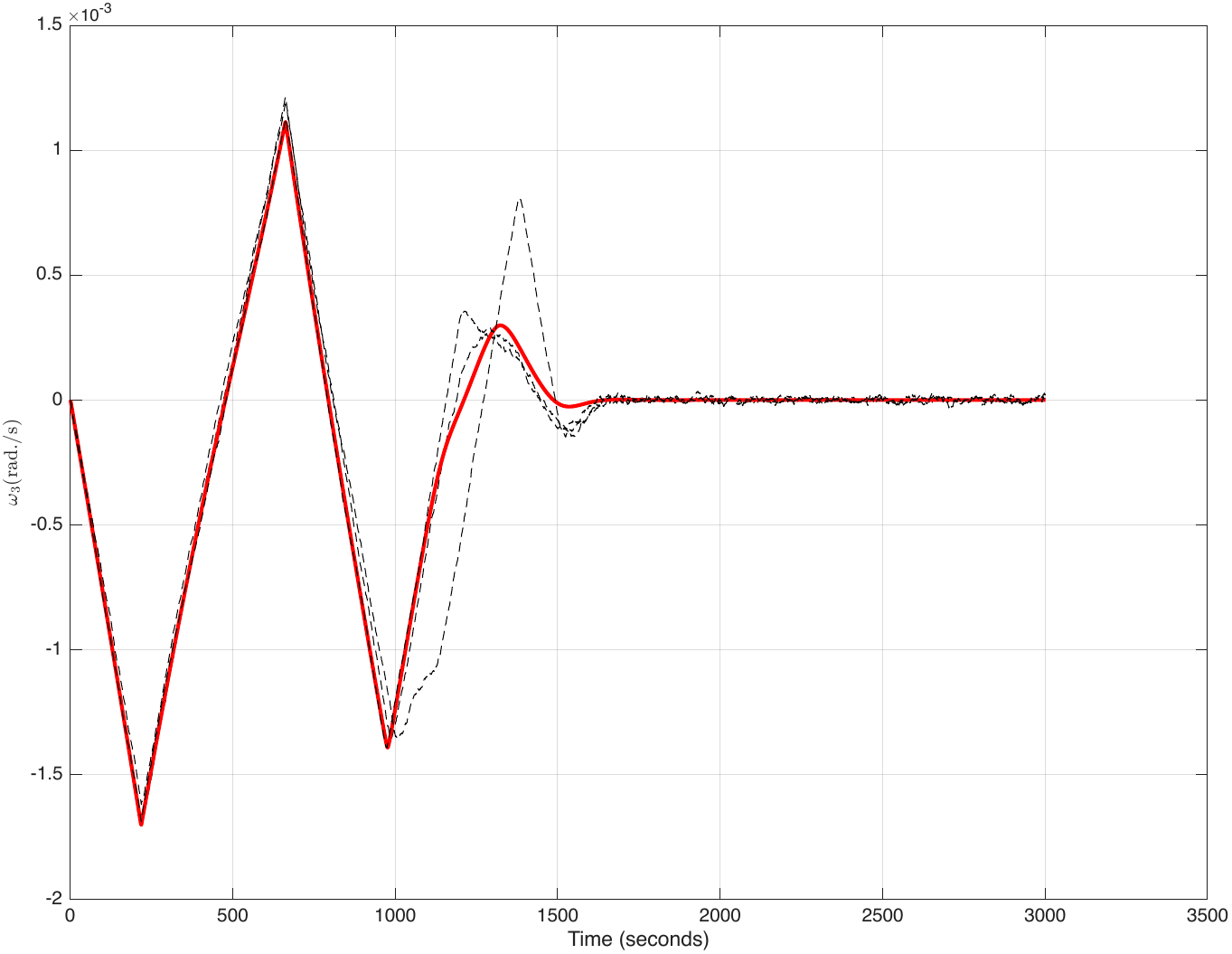}
    \caption{Closed-loop attitude response under disturbances and measurement noise. The nominal trajectory is shown in red and three representative perturbed trajectories are shown using dotted black curves.}
    \label{fig:attitude_closed_loop_perturbed}
\end{figure*}

Without feedback, disturbances accumulate and produce visible deviations from the nominal trajectory in all six states, with larger deviations in the Euler angles than in the angular velocities, as shown in Fig.~\ref{fig:attitude_open_loop_perturbed}. Using the finite-barrier Box-iLQR feedback gains, the perturbed trajectories are driven back toward the nominal solution and terminal errors remain small despite simultaneous process disturbances and measurement noise; three representative realizations are shown in Fig.~\ref{fig:attitude_closed_loop_perturbed}.

For $N_{\mathrm{MC}}=100$ independent realizations, define
\[
\mathbf{e}_{\eta,k}^{(j)}
=
\operatorname{wrap}
\left(
\boldsymbol{\eta}^{(j)}_k-\bar{\boldsymbol{\eta}}_k
\right),
\qquad
\mathbf{e}_{\omega,k}^{(j)}
=
\boldsymbol{\omega}^{(j)}_k-\bar{\boldsymbol{\omega}}_k.
\]

The corresponding mean and maximum errors are
$
e_{z,\mathrm{mean}}=\frac{1}{N_{\mathrm{MC}}}
\sum_{j=1}^{N_{\mathrm{MC}}}
\left|
\mathbf{e}_{z}^{(j)}
\right|_2,
~
e_{z,\max}=\max_j
\left|
\mathbf{e}_{z}^{(j)}
\right|_2,
~ \text{where}~
z\in{\eta,\omega}.
$

\begin{table}[t]
\centering
\caption{Attitude-control errors for 100 disturbance and measurement-noise realizations.}
\label{tab:attitude_mc_results}
\begin{tabular}{lcccc}
\hline
& \multicolumn{2}{c}{Without feedback}
& \multicolumn{2}{c}{With feedback}\\
\cline{2-5}
Error quantity & Mean & Maximum & Mean & Maximum\\
\hline
Euler-angle error [rad]
& $0.0620$ & $0.0702$
& $4.68\times10^{-5}$ & $1.06\times10^{-4}$\\
Angular-rate error [rad/s]
& $3.12\times10^{-5}$ & $3.34\times10^{-5}$
& $1.05\times10^{-5}$ & $2.21\times10^{-5}$\\
\hline
\end{tabular}
\end{table}

Table~\ref{tab:attitude_mc_results} shows the expected trend: Box-iLQR feedback reduces both the mean and worst-case errors relative to open-loop execution.

\subsection{Fuel-Optimal Planar Orbit Transfer}
\label{subsec:orbit_results}

The planar state and controls are

$$
\mathbf{x}=
\begin{bmatrix}
r & u & v
\end{bmatrix}^{T},
\qquad
\mathbf{u}_c=
\begin{bmatrix}
s & \gamma
\end{bmatrix}^{T},
$$

where $r$ is radial position, $u$ radial velocity, $v$ transverse velocity, $s$ normalized thrust magnitude, and $\gamma$ thrust direction. The equations of motion are
\begin{align}
\dot r &=u,\\
\dot u &=\frac{v^2}{r}-\frac{\mu}{r^2}+sa\sin\gamma,\\
\dot v &=-\frac{uv}{r}+sa\cos\gamma,
\end{align}
subject to
$$
0\leq s\leq\bar{s},
\qquad
-\frac{\pi}{3}\leq\gamma\leq\frac{\pi}{3}.
$$

The boundary conditions are

$$
r(0)=1.08,\qquad
u(0)=0,\qquad
v(0)=\sqrt{\frac{\mu}{r(0)}},
$$

and

$$
r(t_f)=1.1,\qquad
u(t_f)=0,\qquad
v(t_f)\;\text{free},\qquad
t_f=6\pi.
$$
The normalized parameters used in the simulation are
$
\mu=1,~~
a=2,~~
\bar{s}=1,
$
with
$
\Delta t=2\pi/10000,~~
N=10000.
$
The feedback gains used for the closed-loop simulations are retained at
$\sigma_{\mathrm{fb}}\approx  10^{-3}$.

The minimum-fuel objective $\int_0^{t_f}s(t),dt$ is approximated using
$
\ell_s(s)=\sqrt{s^2+\alpha^2}-\alpha,
\qquad
\alpha=10^{-12},
$
which approaches $|s|$ away from zero while retaining well-defined first- and second-order derivatives. The discrete objective is
\begin{equation}
\begin{split}
J_h={}
\Delta t
\sum_{k=0}^{N-1}
\left(
\sqrt{s_k^2+\alpha^2}-\alpha
+
\frac{1}{2}\rho_\gamma\gamma_k^2
\right)\
+
\frac{\rho_r}{2}(r_N-1.1)^2
+
\frac{\rho_u}{2}u_N^2,
\end{split}
\label{eq:orbit_penalty_cost}
\end{equation}
where
$
\rho_r=10^9,\qquad
\rho_u=10^6,\qquad
\rho_\gamma=5\times10^{-9}.
$
No terminal penalty is imposed on $v_N$.

\begin{figure*}[t]
\centering
\includegraphics[width=0.45\textwidth]{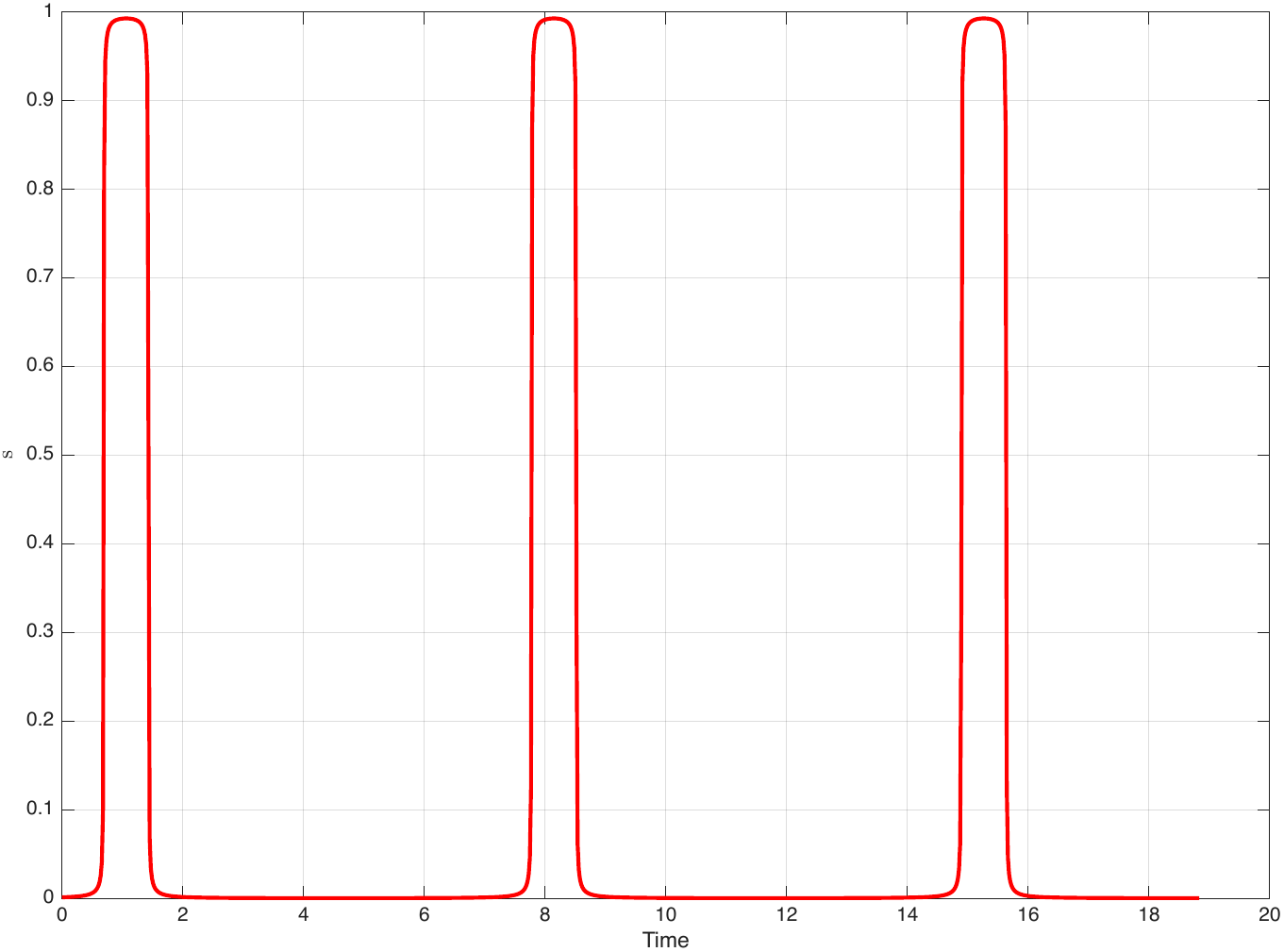}
\hfill
\includegraphics[width=0.45\textwidth]{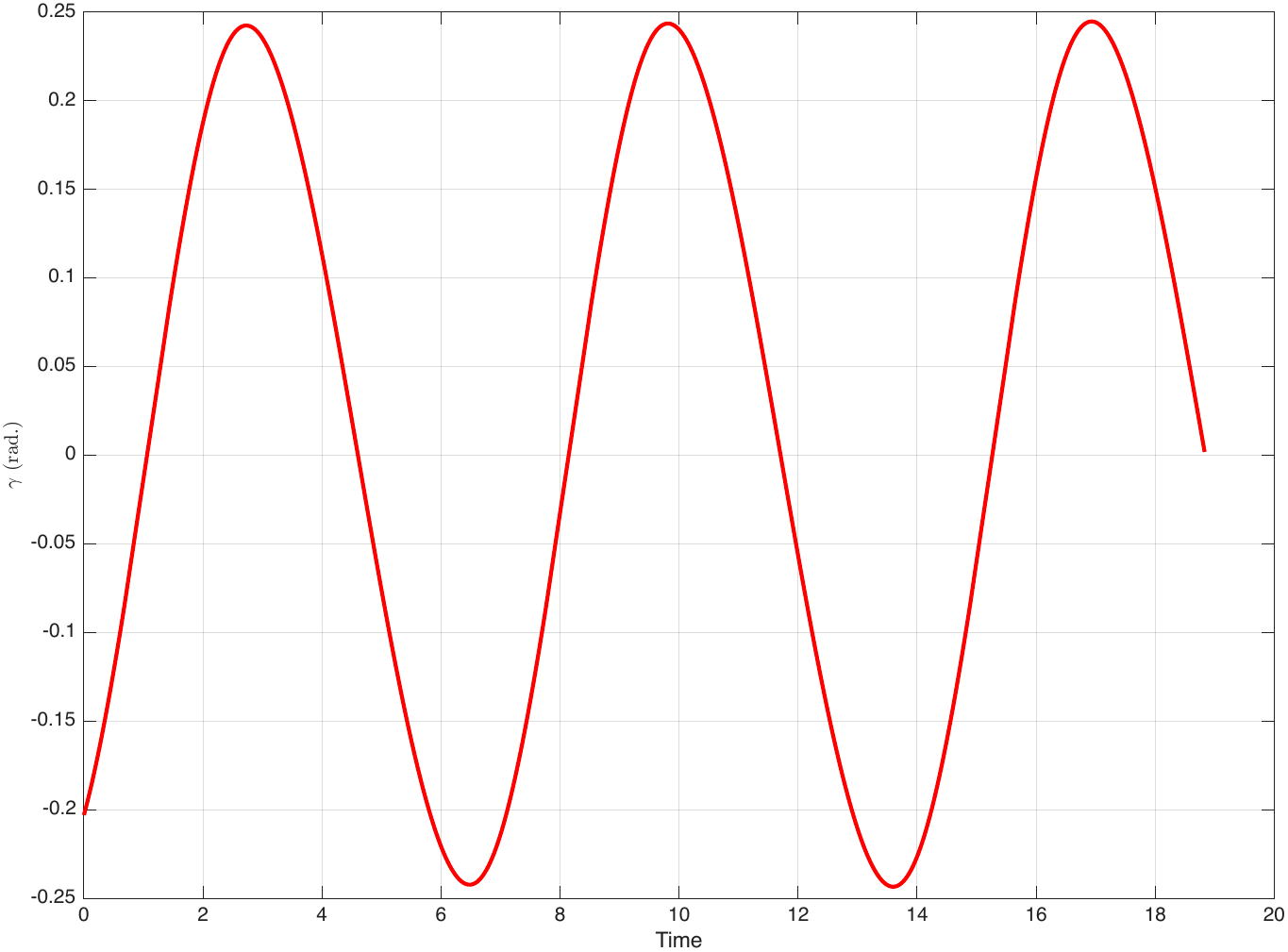}
\caption{Optimal orbit-transfer controls: (a) normalized thrust magnitude and (b) thrust-direction angle $\gamma$.}
\label{fig:orbit_controls}
\end{figure*}

The optimized thrust magnitude in Fig.~\ref{fig:orbit_controls} has a bang--off--bang-type structure with three active-thrust arcs separated by coast intervals; each principal thrust arc lasts approximately $0.7$ normalized time units. The thrust direction repeatedly approaches the bounds $\gamma=\pm\pi/3$. 


\begin{figure}[t]
\centering
\includegraphics[width=0.44\textwidth]{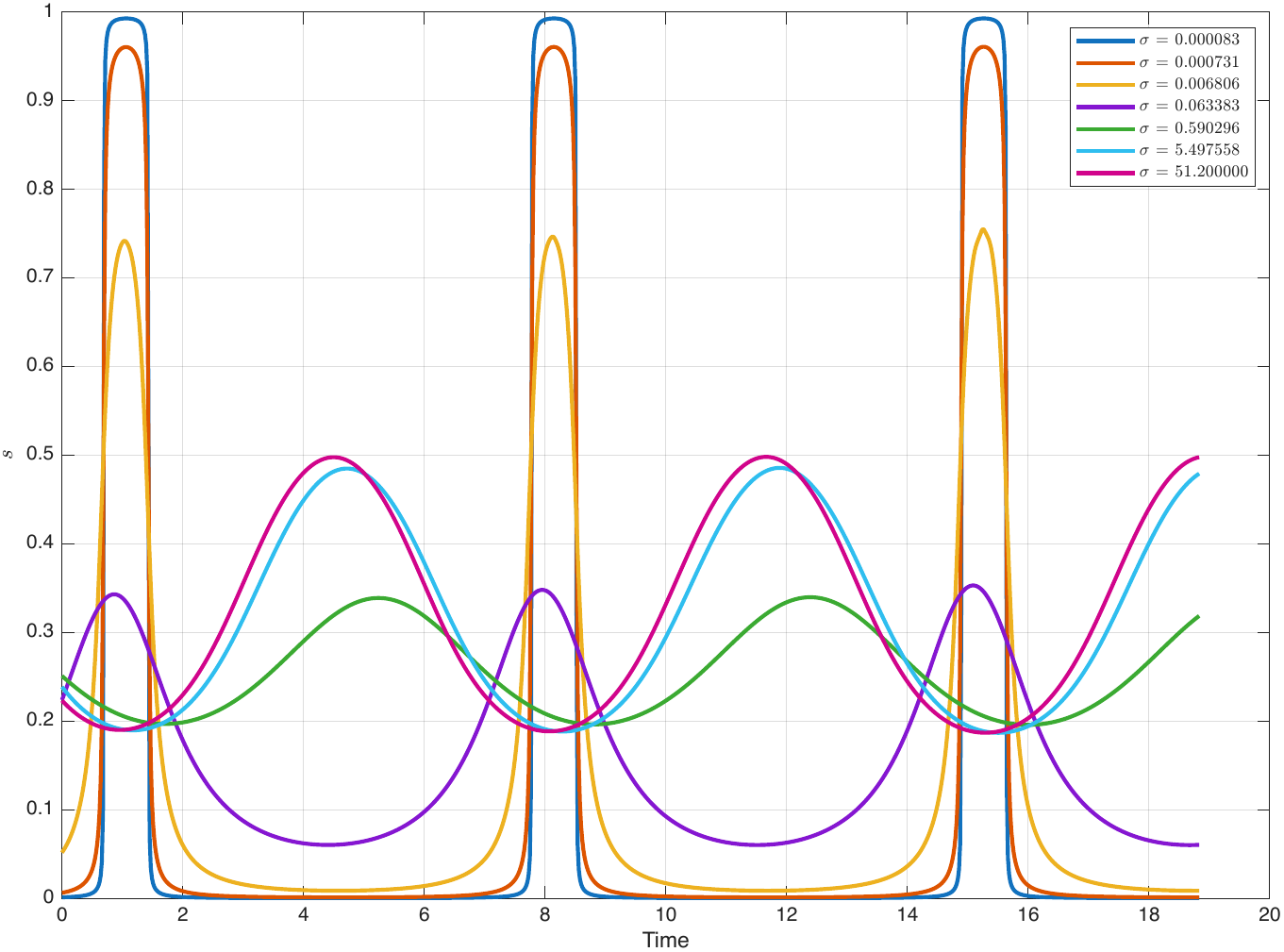}
\hfill
\includegraphics[width=0.44\textwidth]{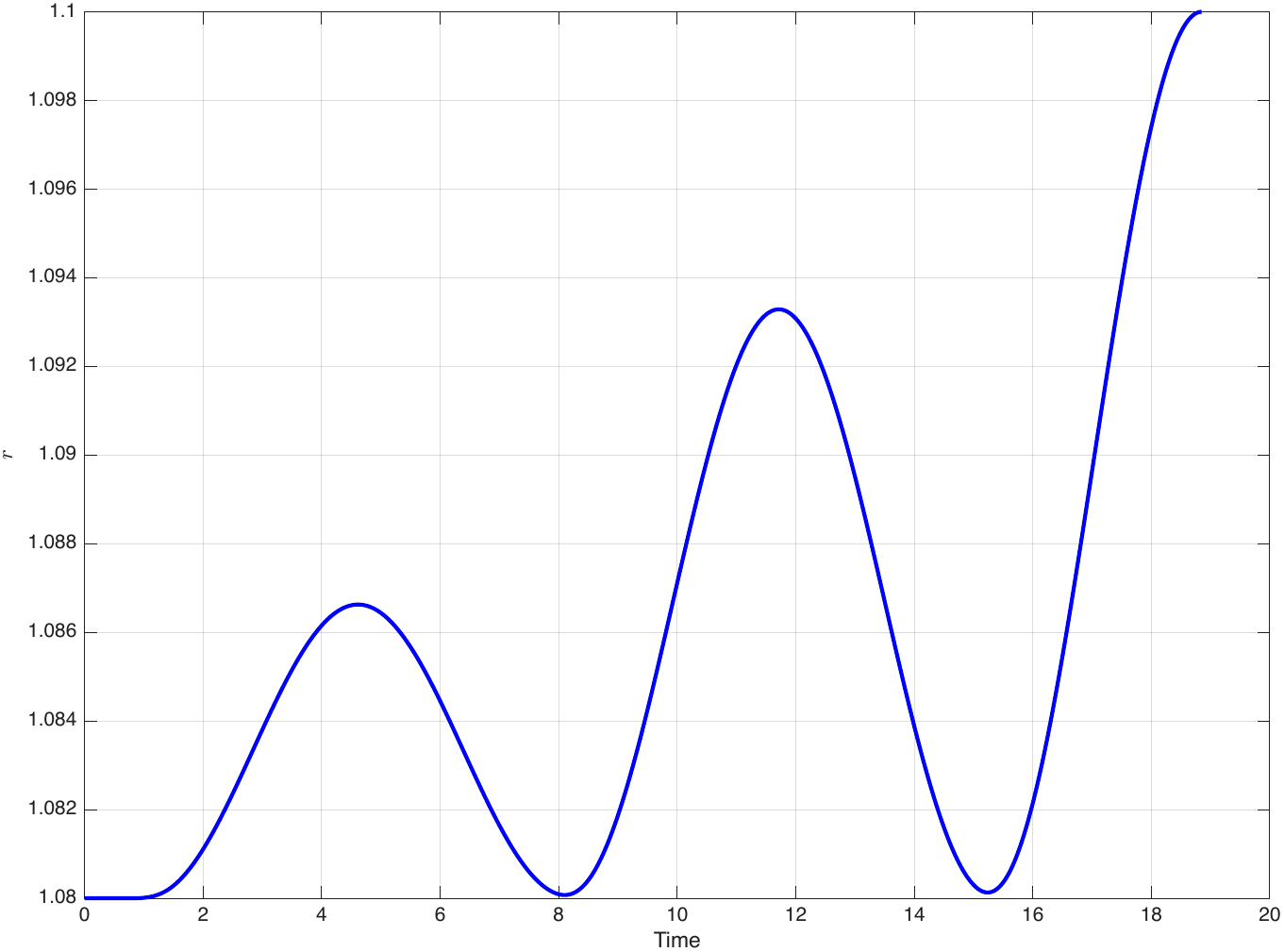}
\includegraphics[width=0.44\textwidth]{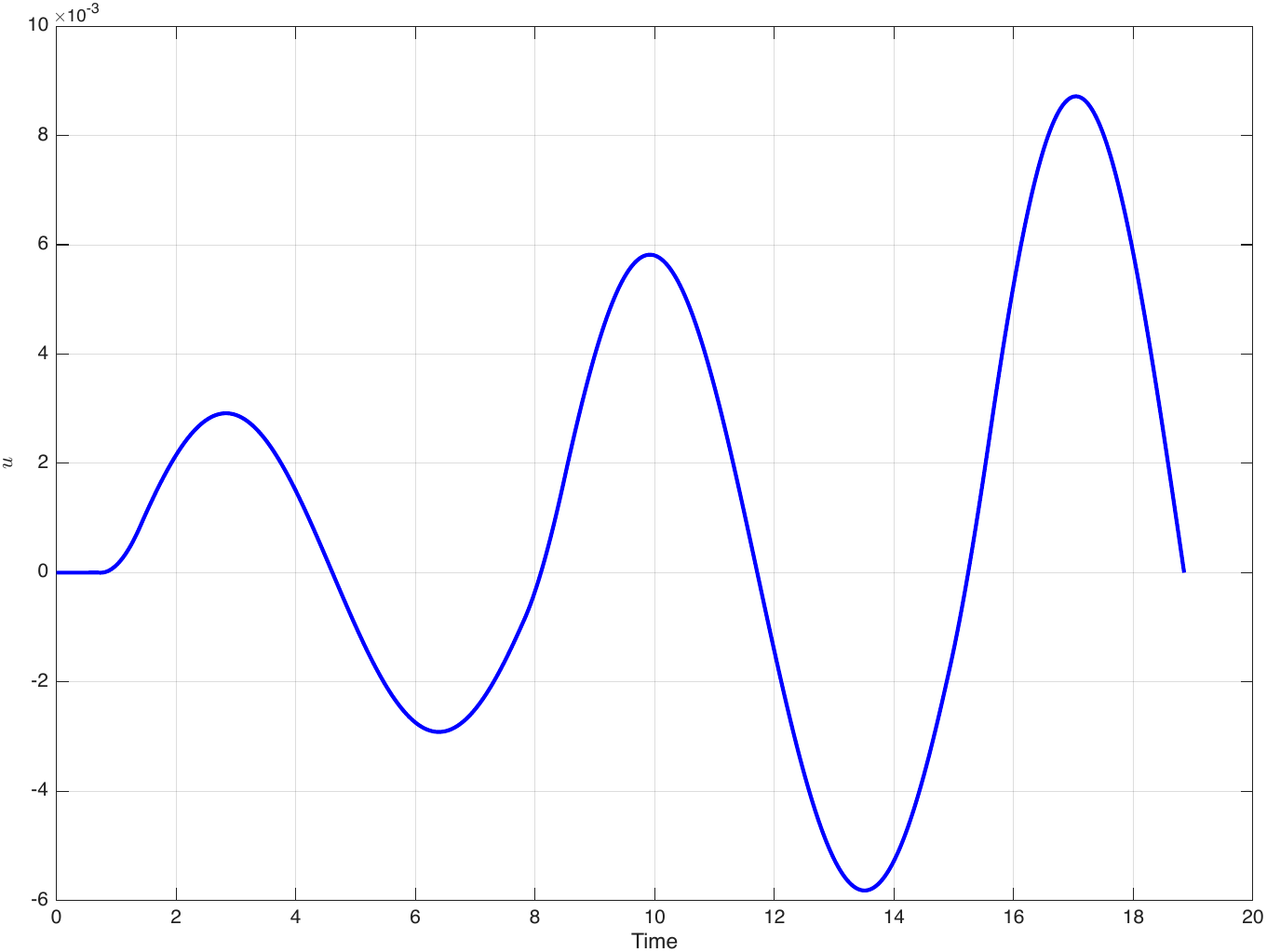}
\hfill
\includegraphics[width=0.44\textwidth]{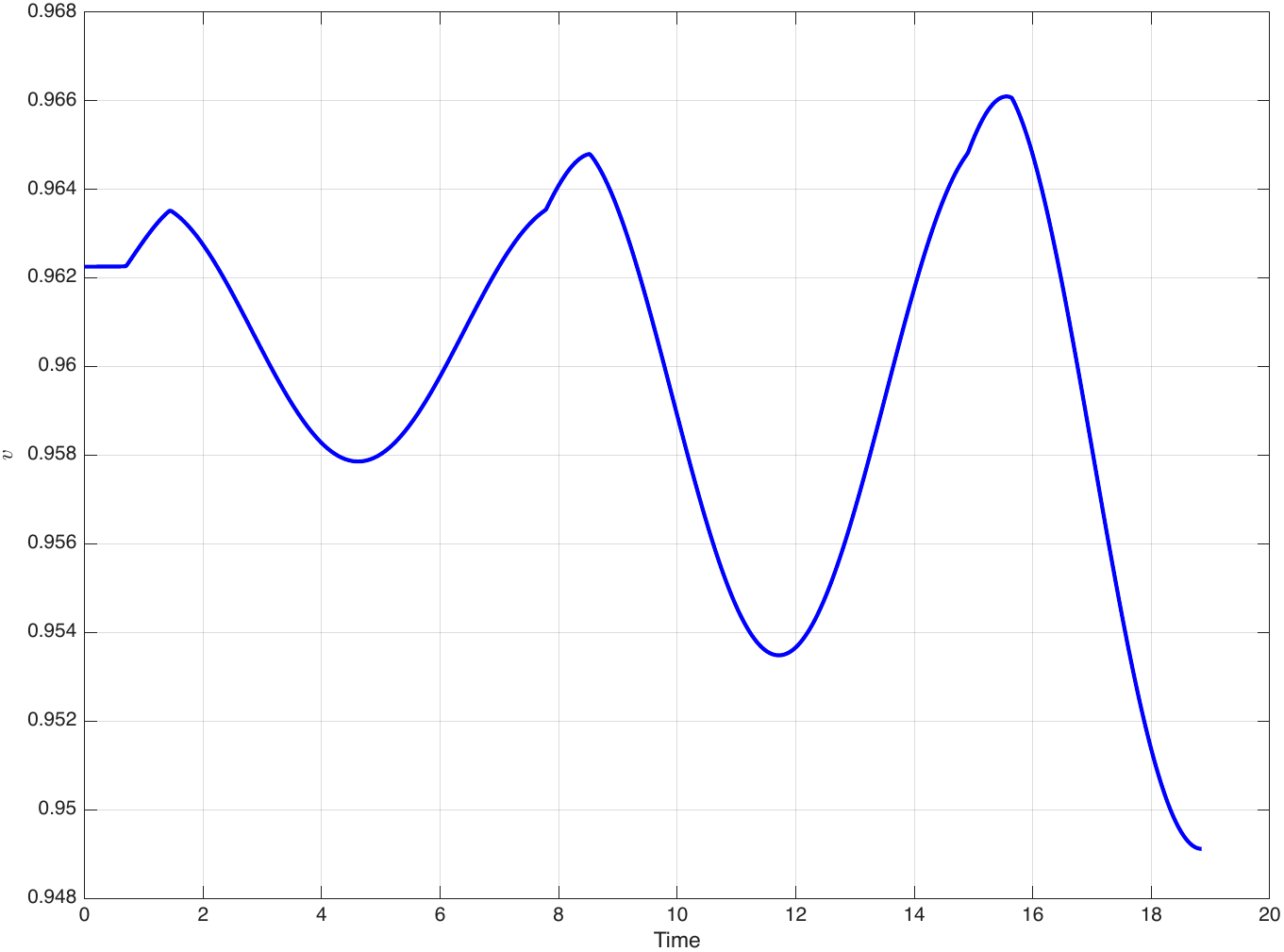}
\caption{Evolution of the thrust-magnitude history for decreasing values of the logarithmic-barrier parameter, together with the corresponding nominal orbit-transfer state histories: radial position $r$, radial velocity $u$, and transverse velocity $v$.}
\label{fig:orbit_states}
\end{figure}

Figure~\ref{fig:orbit_states} shows the evolution of the thrust-magnitude profile as the logarithmic-barrier parameter is reduced, together with the corresponding nominal state histories. As the barrier parameter decreases, the thrust profile approaches the limiting bang--off--bang structure. The radial position completes approximately $2.5$ oscillations with increasing amplitude while transitioning from $r(0)=1.08$ to $r_f=1.1$, while the radial velocity initially grows in amplitude before decaying to zero. The terminal residuals are

$$
r(t_f)-r_f\approx-0.1115\times10^{-6},
\qquad
u(t_f)-u_f\approx-0.0429\times10^{-6}.
$$

Thus, the quadratic terminal penalties accurately approximate the terminal equality constraints while remaining compatible with the local quadratic expansions used by iLQR.

Control uncertainty is introduced as
\begin{align*}
\widetilde{\gamma}_k
=
\gamma_k+\delta\gamma_k,\quad
\delta\gamma_k
\sim
\mathcal{N}(0,\sigma_\gamma^2),\quad
\sigma_\gamma
=3^\circ,\quad
\widetilde{s}_k
=
s_k+\delta s_k,\quad
\delta s_k
\sim
\mathcal{N}\left(0,(0.1s_k)^2\right).
\end{align*}
All perturbed controls are projected onto their admissible bounds before propagation.

The robustness results are summarized in Table~\ref{tab:orbit_mc_results}. Without feedback, both terminal radial-position and radial-velocity errors are of order $10^{-4}$, several orders of magnitude larger than the nominal residuals. With Box-iLQR feedback, the radial-position error decreases to order $10^{-7}$, corresponding under the adopted normalization to dimensional errors of order tens of meters, while the radial-velocity error decreases to approximately order $10^{-6}$.

\begin{table}[t]
    \centering
    \caption{Planar orbit-transfer robustness under control uncertainty.}
    \label{tab:orbit_mc_results}
    \begin{tabular}{lcc}
        \hline
        Case & Radial-position error & Radial-velocity error\\
        \hline
        Nominal
        & $-0.1115\times10^{-6}$
        & $-0.0429\times10^{-6}$\\
        No feedback
        & $\mathcal{O}(10^{-4})$
        & $\mathcal{O}(10^{-4})$\\
        Box-iLQR feedback
        & $\mathcal{O}(10^{-6})$
        & $\mathcal{O}(10^{-6})$\\
        \hline
    \end{tabular}
\end{table}

The actual fuel-consumption metric
$
\int_0^{t_f}s(t)\,dt
$
is distinct from the smoothed cost used internally by Box-iLQR. The nominal optimal fuel cost is $2.2228$, and for every independent perturbation realization the feedback-corrected fuel cost remains less than $1\%$ above this value. Thus, the feedback policy substantially reduces terminal error while introducing only a small increase in actual fuel consumption.

\subsection{Fixed-Time $L_2$-to-$L_2$ Halo-Orbit Transfer}
\label{subsec:cr3bp_halo_transfer}

The third example considers a transfer between two southern $L_2$ Halo orbits in the Earth--Moon system. The boundary conditions, transfer time, and reference trajectory follow Refs.~\cite{Saloglu2024AccelerationBased,Yamamoto2025JacobianColoring,Yamamoto2026NonlinearMPC}. Unlike the impulsive and trajectory-correction formulations used in those works, the present maneuver is solved directly as a continuous control-constrained transfer. The fixed final time is
$
t_f=13.0789122308~\mathrm{TU},
$
and free-final-time optimization is not considered.

The CR3BP state and control are
$
\mathbf{x}
=
\begin{bmatrix}
x&y&z&\dot{x}&\dot{y}&\dot{z}
\end{bmatrix}^{T},
~~
\mathbf{u}
=
\begin{bmatrix}
u_x&u_y&u_z
\end{bmatrix}^{T},
$
with position in LU, velocity in VU, time in TU, and acceleration in $\mathrm{LU}/\mathrm{TU}^2$. The controlled equations are
\begin{align}
\ddot{x}-2\dot{y}
&=
\frac{\partial\Omega}{\partial x}+u_x,\\
\ddot{y}+2\dot{x}
&=
\frac{\partial\Omega}{\partial y}+u_y,\\
\ddot{z}
&=
\frac{\partial\Omega}{\partial z}+u_z,
\end{align}
where

$$
\Omega(x,y,z)
=
\frac{1-\mu}{d_1}
+
\frac{\mu}{d_2}
+
\frac{1}{2}(x^2+y^2),
$$

$$
d_1=\sqrt{(x+\mu)^2+y^2+z^2},
\qquad
d_2=\sqrt{(x-1+\mu)^2+y^2+z^2}.
$$

Each control component satisfies
\begin{equation}
-\frac{0.05}{\sqrt{3}}
\leq
u_i(t)
\leq
\frac{0.05}{\sqrt{3}},
\qquad
i\in{x,y,z}.
\label{eq:cr3bp_control_bounds}
\end{equation}

A smooth minimum-$\Delta v$ running cost is used,
\begin{equation}
\ell(\mathbf{u})=100
\left(
\sqrt{u_x^2+u_y^2+u_z^2+\alpha^2}
-\alpha
\right),
\qquad
\alpha=10^{-6},
\label{eq:cr3bp_smooth_fuel_cost}
\end{equation}
which approaches $100|\mathbf{u}|$ for $|\mathbf{u}|\gg\alpha$ while remaining twice differentiable near zero. The objective is
\begin{equation}
J=
\frac{1}{2}
(\mathbf{x}(t_f)-\mathbf{x}_f)^T
Q_f
(\mathbf{x}(t_f)-\mathbf{x}_f)
+
\int_0^{t_f}\ell(\mathbf{u}),dt,
\end{equation}
with
\begin{equation}
Q_f=50\times10^4
\operatorname{diag}
(10^4,10^4,10^4,1,1,1).
\end{equation}
Because the final state is imposed using a quadratic penalty rather than an equality constraint, a small terminal residual remains.

\begin{figure}[t]
\centering
\includegraphics[width=0.6\columnwidth]{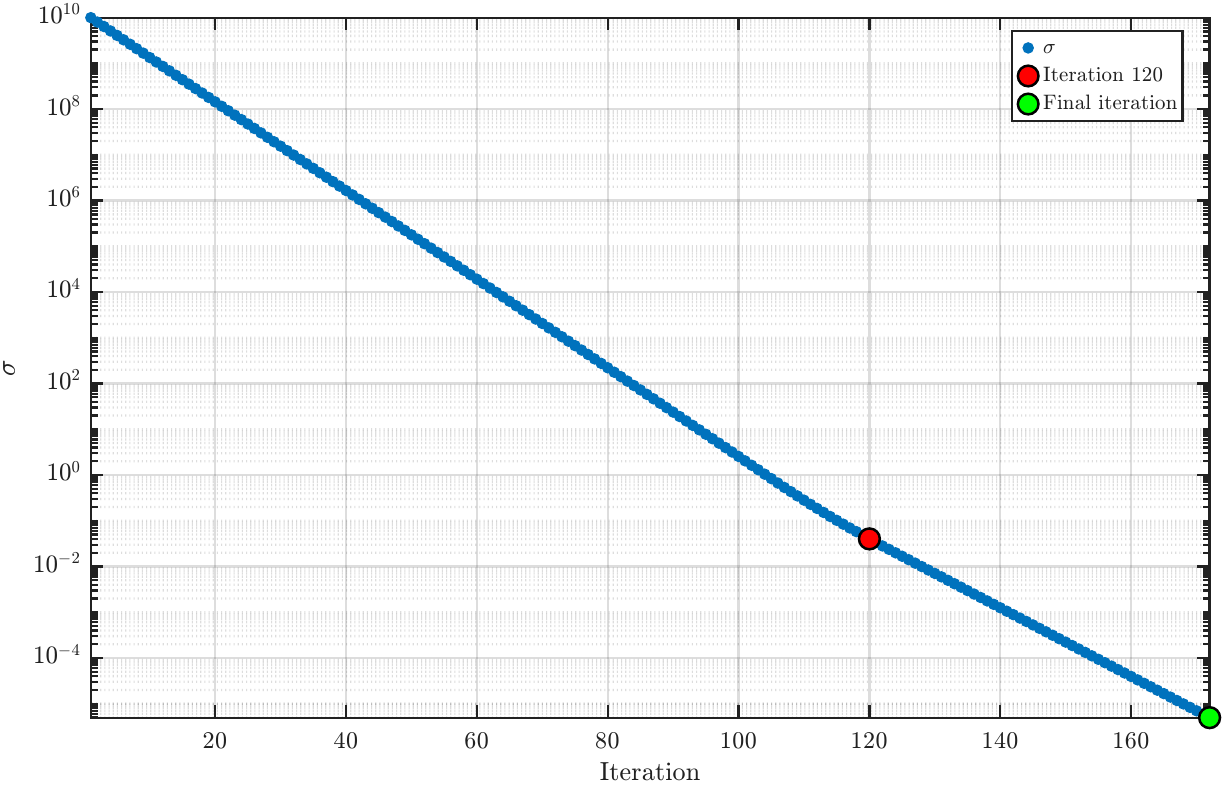}
\caption{Reduction of the control log-barrier parameter. Iteration 120 corresponds to $\sigma=0.0401$ and iteration 172 to $\sigma=4.945\times10^{-6}$.}
\label{fig:cr3bp_barrier}
\end{figure}

Figure~\ref{fig:cr3bp_barrier} shows the control-barrier continuation over 172 outer iterations. Two solutions are considered: iteration 120 with
$
\sigma=0.0401
$
and the final iteration with
$
\sigma=4.945\times10^{-6}.
$

The final value more closely approximates the fixed-time fuel-optimal constrained solution, whereas $\sigma=0.0401$ retains greater interior control authority and produces substantially more effective feedback gains.

At the final barrier value, all three controls contain extended coast intervals; $u_y$ and $u_z$ additionally contain saturated arcs, while $u_x$ remains strictly interior. After removing the factor of 100 in Eq.~\eqref{eq:cr3bp_smooth_fuel_cost} and converting LU/TU to dimensional velocity,

$$
\Delta v_{\mathrm{final}}=283.7073~\mathrm{m/s}.
$$

\begin{figure*}[t]
\centering
\subfloat[$u_x$]{
\includegraphics[width=0.31\textwidth]{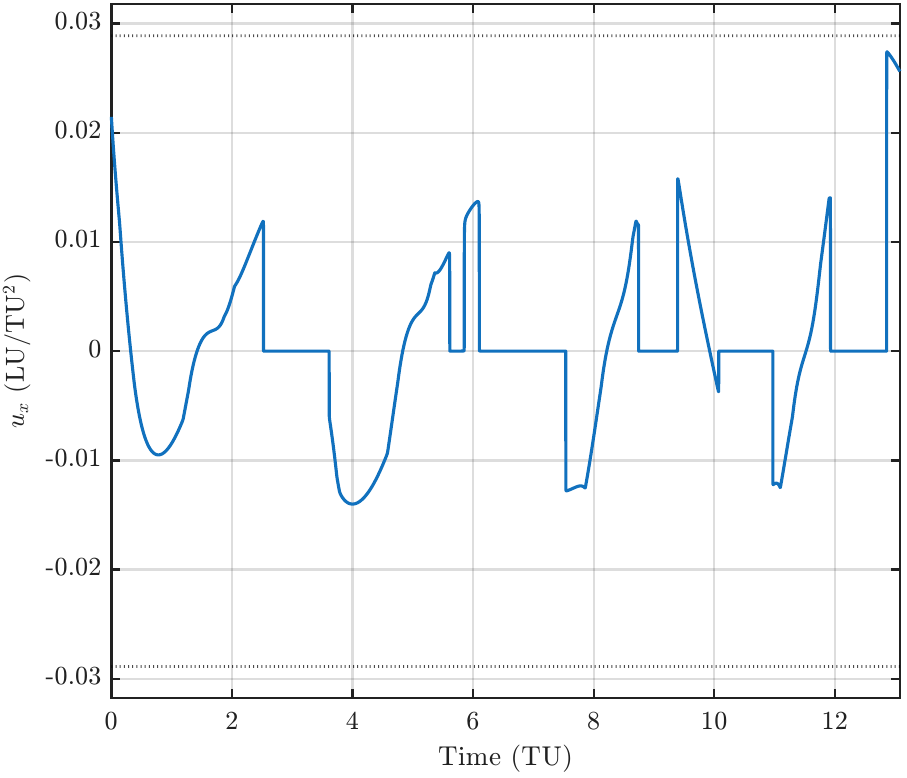}}
\hfill
\subfloat[$u_y$]{
\includegraphics[width=0.31\textwidth]{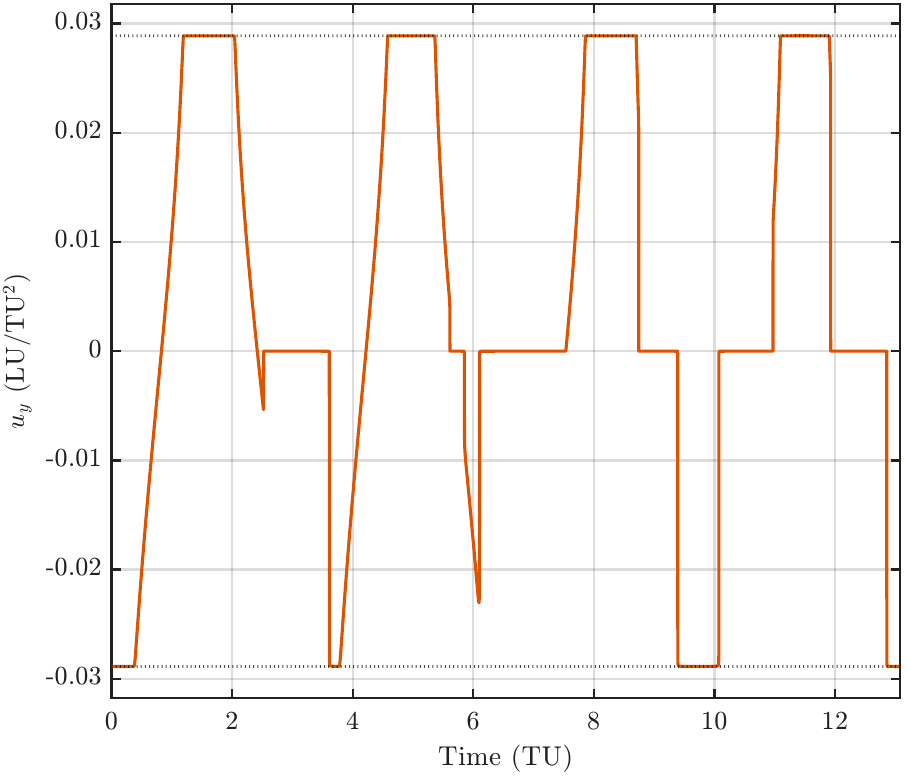}}
\hfill
\subfloat[$u_z$]{
\includegraphics[width=0.31\textwidth]{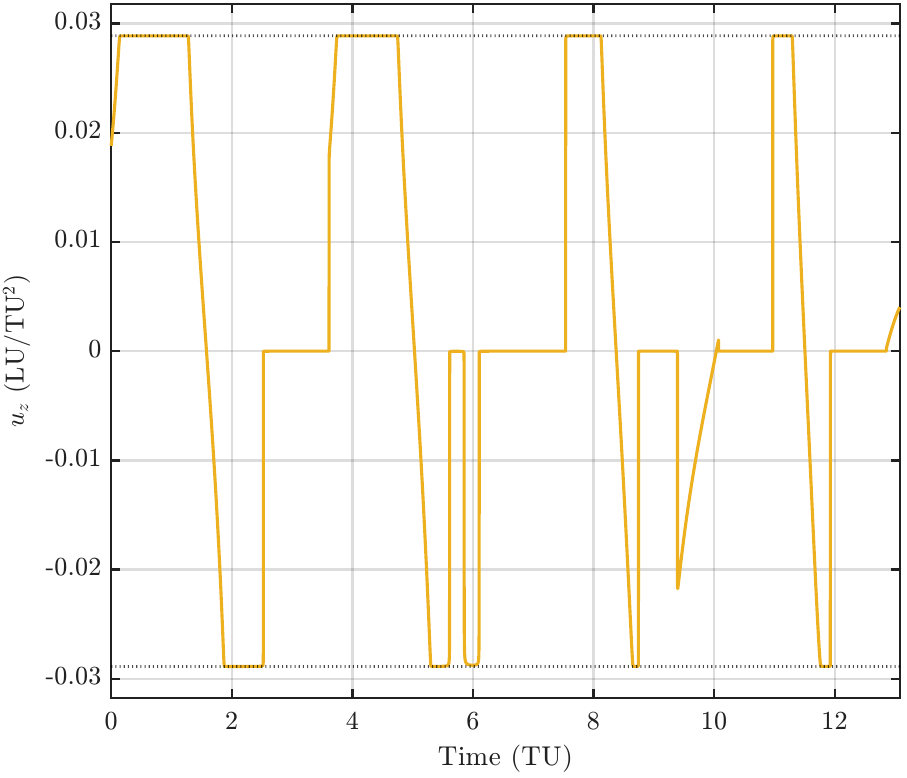}}
\caption{Control histories at $\sigma=4.945\times10^{-6}$. Dotted black lines denote the bounds $\pm0.05/\sqrt{3}~\mathrm{LU}/\mathrm{TU}^2$. The $u_y$ and $u_z$ controls contain saturated arcs and all three controls contain coast intervals.}
\label{fig:cr3bp_final_controls}
\end{figure*}

At iteration 120, the controls still contain coast intervals but remain strictly within the bounds. The maneuver cost is

$$
\Delta v_{120}=287.1338~\mathrm{m/s},
$$

only approximately $1.21\%$ above the final low-barrier value. Thus, a small loss in nominal fuel optimality preserves appreciable control authority for feedback correction.

\begin{figure*}[t]
\centering
\subfloat[$u_x$]{
\includegraphics[width=0.31\textwidth]{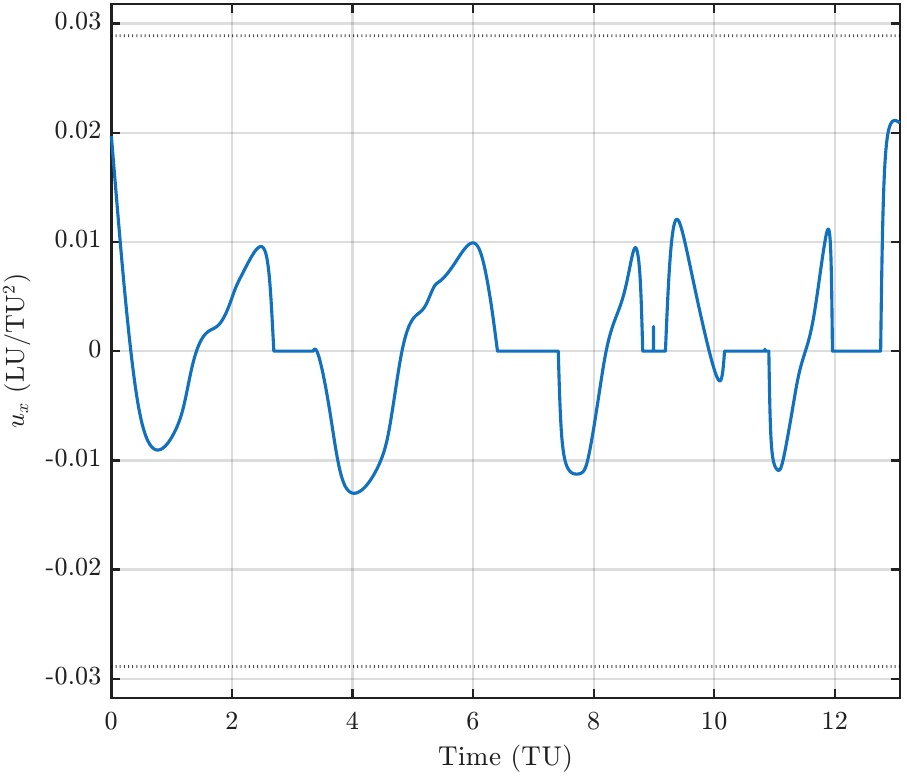}}
\hfill
\subfloat[$u_y$]{
\includegraphics[width=0.31\textwidth]{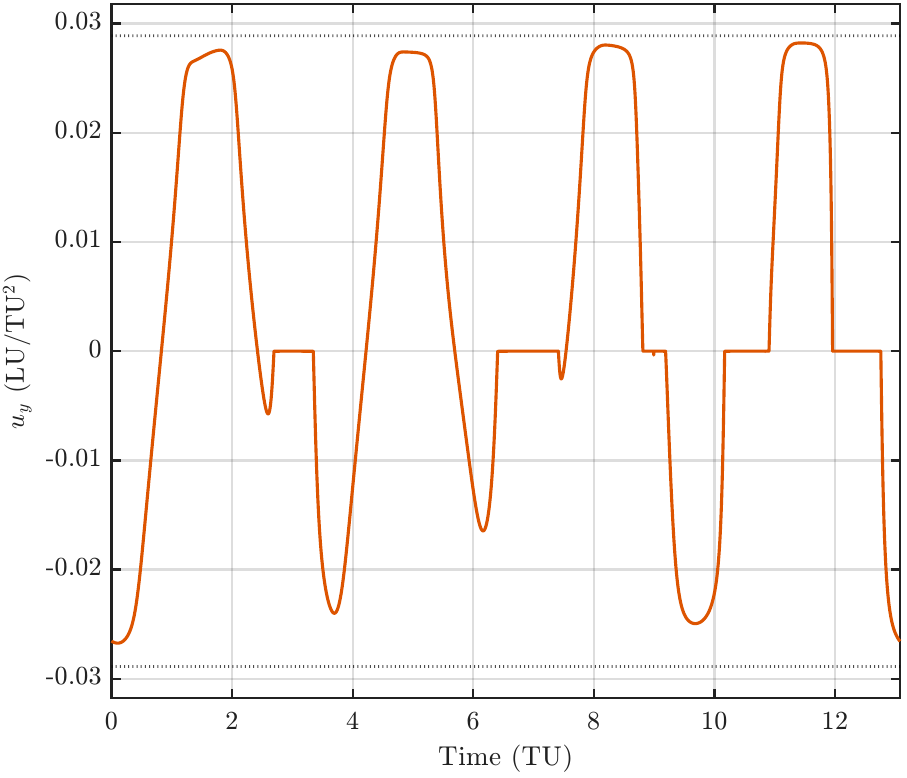}}
\hfill
\subfloat[$u_z$]{
\includegraphics[width=0.31\textwidth]{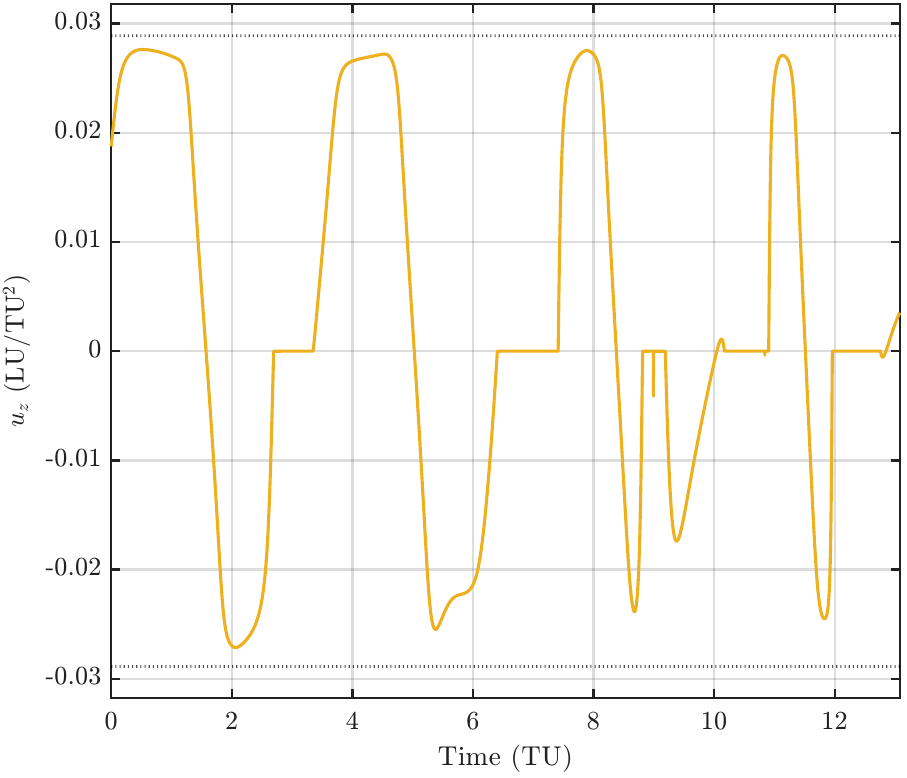}}
\caption{Control histories at iteration 120 ($\sigma=0.0401$). The controls contain coast intervals but remain strictly interior to their bounds.}
\label{fig:cr3bp_120_controls}
\end{figure*}

\begin{figure*}[t]
\centering
\subfloat[$x$]{
\includegraphics[width=0.31\textwidth]{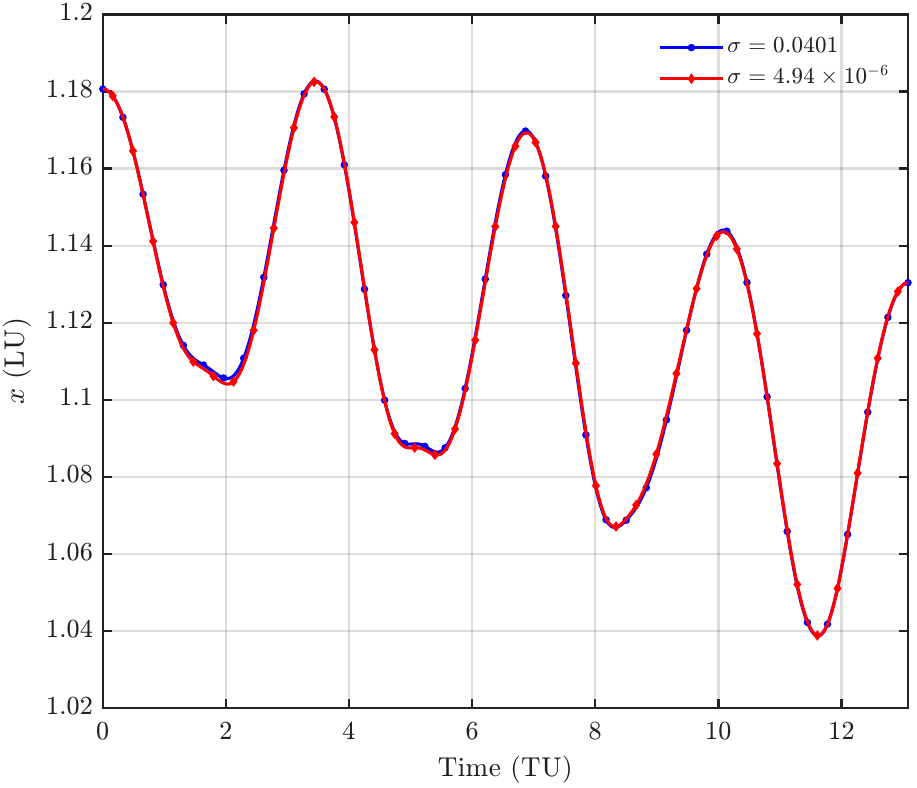}}
\hfill
\subfloat[$y$]{
\includegraphics[width=0.31\textwidth]{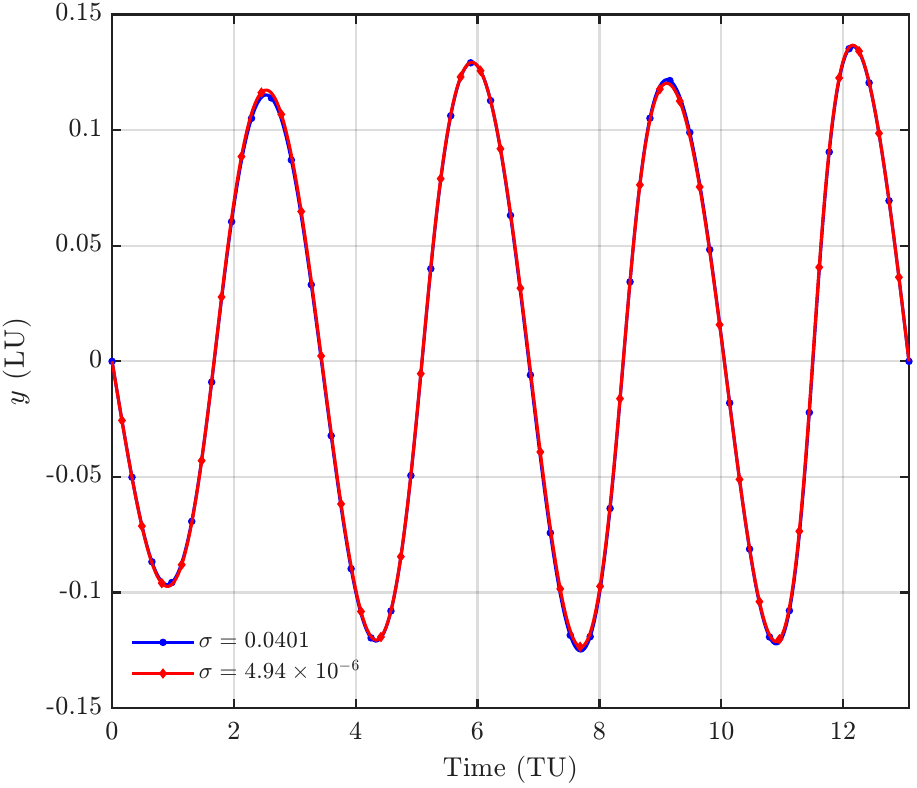}}
\hfill
\subfloat[$z$]{
\includegraphics[width=0.31\textwidth]{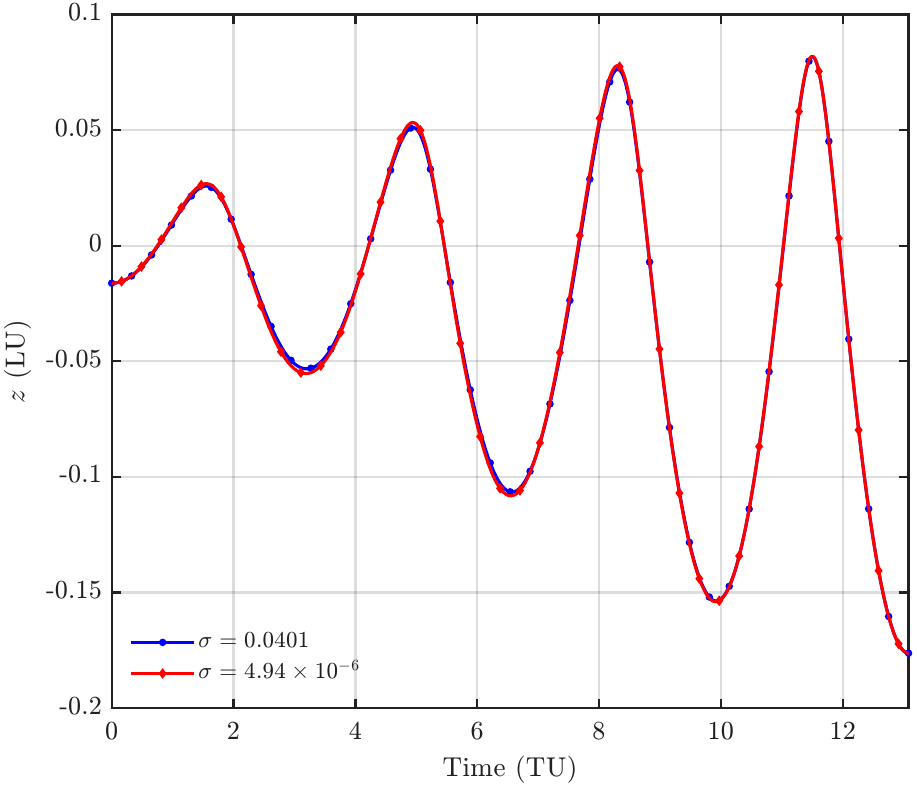}}
\hfill
\subfloat[$\dot{x}$]{
\includegraphics[width=0.31\textwidth]{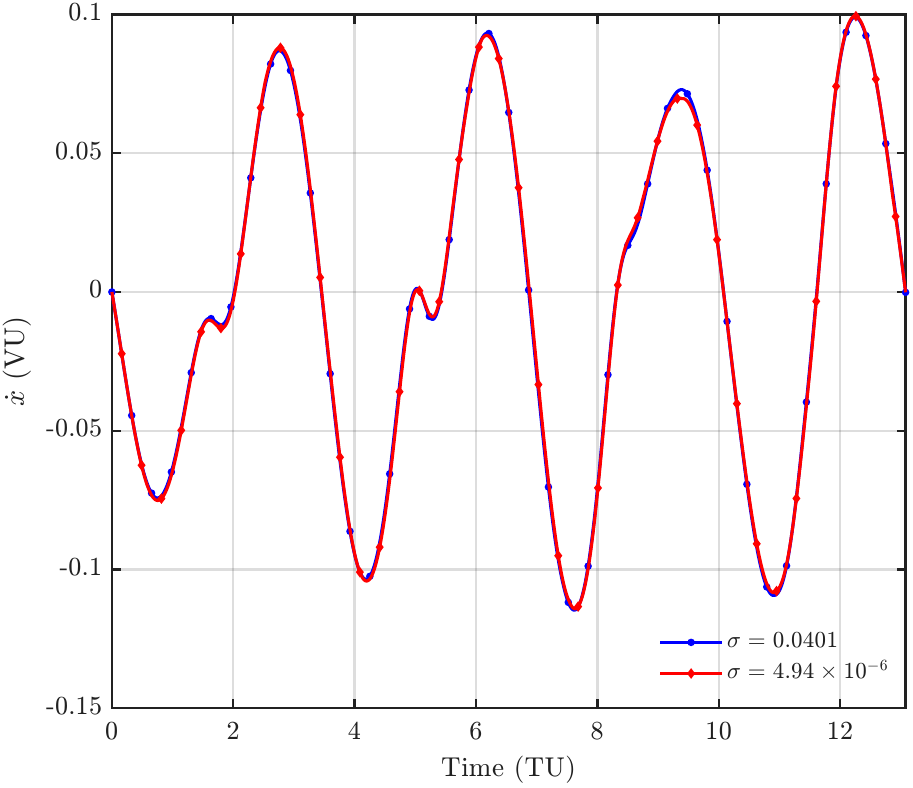}}
\hfill
\subfloat[$\dot{y}$]{
\includegraphics[width=0.31\textwidth]{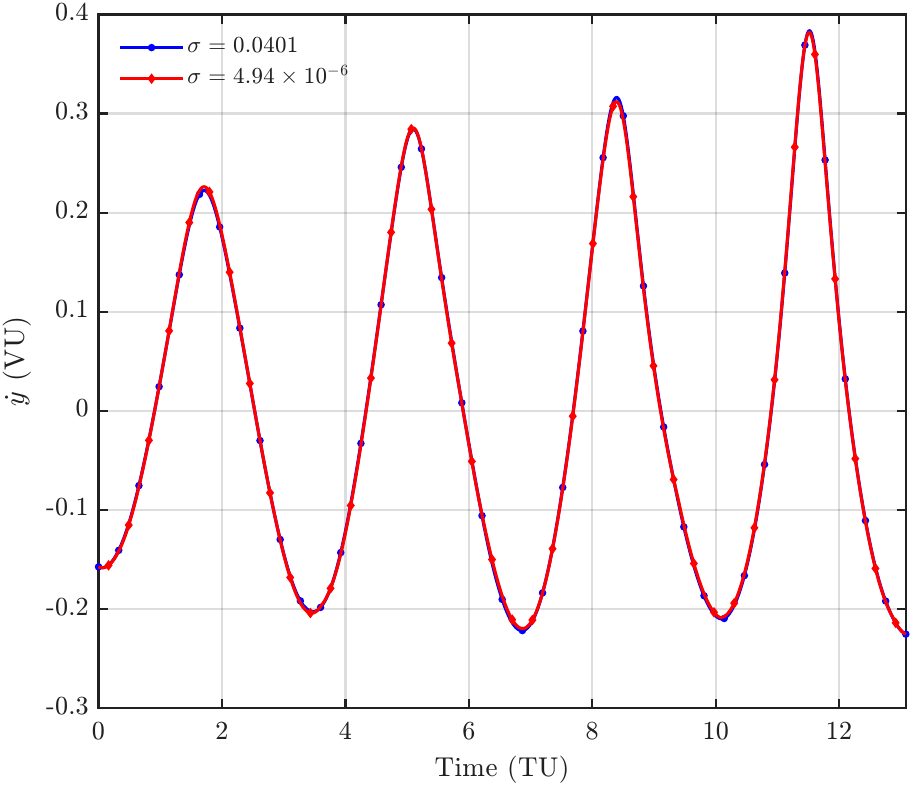}}
\hfill
\subfloat[$\dot{z}$]{
\includegraphics[width=0.31\textwidth]{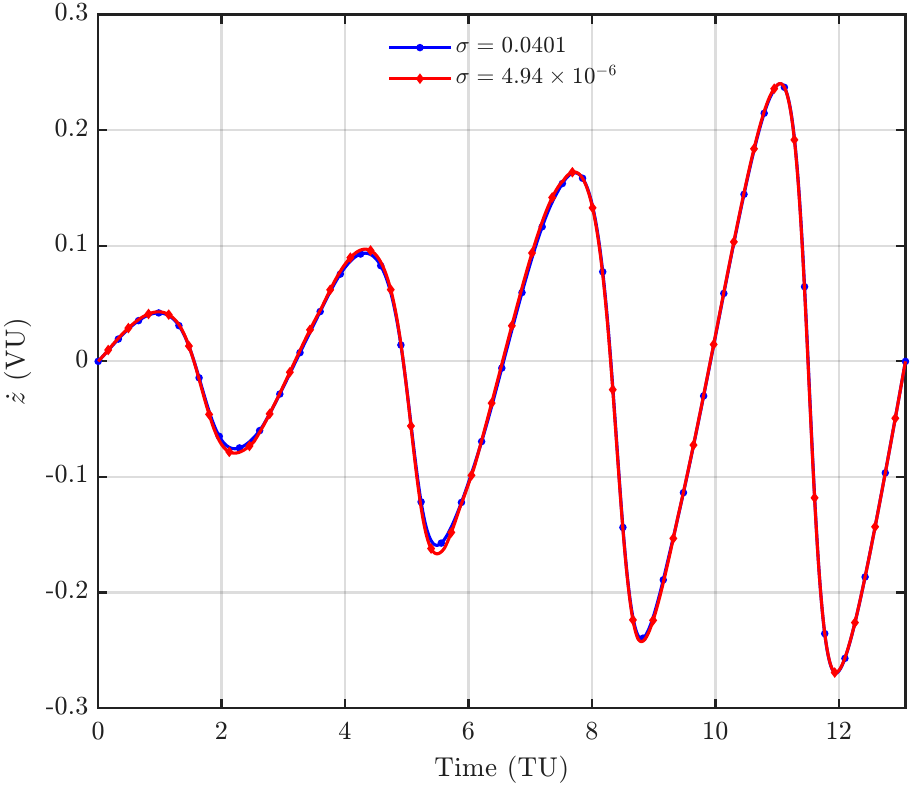}}
\caption{Comparison of nominal CR3BP state histories for $\sigma=0.0401$ and $\sigma=4.945\times10^{-6}$. The former is shown using blue lines with circular markers and the latter using red lines with diamond markers. Position is expressed in LU and velocity in VU.}
\label{fig:cr3bp_state_comparison}
\end{figure*}

Despite their different control margins and closed-loop behavior, the nominal state trajectories at the two barrier values are nearly identical, as shown in Fig.~\ref{fig:cr3bp_state_comparison}; all six position and velocity histories overlap over nearly the entire transfer.

The dimensional terminal residual at $\sigma=0.0401$ is
\begin{equation}
\mathbf{e}_{f,120}=\begin{bmatrix}
13.3345&
4.5538&
2.6976&
-0.1365&
0.2022&
-0.0222
\end{bmatrix}^{T},
\end{equation}
where the first three components are in meters and the final three in meters per second, giving

$$
\|\mathbf{e}_{r,f,120}\|=14.3465~\mathrm{m},
\qquad
\|\mathbf{e}_{v,f,120}\|=0.2450~\mathrm{m/s}.
$$

At $\sigma=4.945\times10^{-6}$,
\begin{equation}
\mathbf{e}_{f,\mathrm{final}}=\begin{bmatrix}
12.7258&
4.4417&
2.4579&
-0.1353&
0.1866&
-0.0213
\end{bmatrix}^{T},
\end{equation}
with

$$
\|\mathbf{e}_{r,f,\mathrm{final}}\|=13.7009~\mathrm{m},
\qquad
\|\mathbf{e}_{v,f,\mathrm{final}}\|=0.2315~\mathrm{m/s}.
$$

The slightly smaller final-barrier residuals are consistent with its lower nominal objective, although the difference is small relative to the difference in closed-loop robustness.

A total of 100 independent Monte Carlo realizations are performed. Each control component is perturbed according to
\begin{equation}
\widetilde{u}_{i,k}
\sim
\mathcal{N}
\left(
\bar{u}_{i,k},
\left(
\frac{|\bar{u}_{i,k}|}{30}
\right)^2
\right),
\qquad i\in{x,y,z},
\label{eq:cr3bp_control_noise}
\end{equation}
so that the $3\sigma$ uncertainty equals $10\%$ of the magnitude of the corresponding nominal control. The noisy control is corrected using gains obtained from the same Box-iLQR iteration,
\begin{equation}
\mathbf{u}^{\mathrm{cl}}_k=\Pi_{\mathcal U}
\left[
\widetilde{\mathbf{u}}_k
+
\mathbf{K}_k
\left(
\mathbf{x}^{\mathrm{noisy}}_k-\bar{\mathbf{x}}_k
\right)
\right].
\end{equation}

Table~\ref{tab:cr3bp_terminal_errors} summarizes the dimensional terminal errors. Without feedback, control perturbations cause complete departure from the nominal transfer, with position errors of order $10^9$ m. Feedback from the final low-barrier solution also performs poorly, with position errors of order $10^8$ m and velocity errors of order $10^3$ m/s. In contrast, the gains retained at $\sigma=0.0401$ reduce the mean and maximum position errors to $20.1178$ m and $55.4886$ m and the corresponding velocity errors to $0.2457$ m/s and $0.2523$ m/s.

\begin{table*}[t]
\centering
\caption{Terminal-error norms for 100 noisy-control realizations.}
\label{tab:cr3bp_terminal_errors}
\renewcommand{\arraystretch}{1.2}
\small
\begin{tabular}{@{}lcccc@{}}
\hline
Case
& \shortstack{Mean position\\error (m)}
& \shortstack{Maximum position\\error (m)}
& \shortstack{Mean velocity\\error (m/s)}
& \shortstack{Maximum velocity\\error (m/s)}\\
\hline
No feedback, $\sigma=0.0401$
& $1.0840\times10^{9}$
& $1.7167\times10^{9}$
& $8.9771\times10^{3}$
& $6.6314\times10^{5}$\\

Feedback, $\sigma=0.0401$
& $20.1178$
& $55.4886$
& $0.2457$
& $0.2523$\\

Feedback, $\sigma=4.945\times10^{-6}$
& $4.9669\times10^{8}$
& $5.5826\times10^{8}$
& $1.5214\times10^{3}$
& $1.7420\times10^{3}$\\
\hline
\end{tabular}
\end{table*}
For the 100 feedback-corrected realizations at $\sigma=0.0401$, the maximum maneuver cost is
$
\Delta v_{\max}=287.3046~\mathrm{m/s},
$
only $0.0595\%$ above the nominal iteration-120 value and approximately $1.27\%$ above the final low-barrier value. The feedback correction therefore adds negligible fuel cost while substantially improving terminal accuracy.

For comparison, Yamamoto et al.~\cite{Yamamoto2026NonlinearMPC} apply NMPC corrections over individual transfer segments and include an additional orbital period after the final impulse to capture the target Halo orbit. Their Monte Carlo study reports a maximum terminal miss distance of approximately $0.6101~\mathrm{km}$, median and maximum correction costs of $308.54~\mathrm{m/s}$ and $409.93~\mathrm{m/s}$, respectively, and a reference impulsive cost of $202.2024~\mathrm{m/s}$. Their maximum corrected cost is therefore approximately twice the reference cost. In the present study, the maximum position error is $55.4886~\mathrm{m}$ without an additional target-orbit correction period, while the increase in $\Delta v$ relative to the selected nominal trajectory is negligible. Because the uncertainty models and correction architectures differ, this comparison is not one-to-one, but it demonstrates the closed-loop performance obtained directly from the Box-iLQR feedback policy.

The CR3BP example also illustrates the distinction between nominal optimality and closed-loop robustness. Reducing the barrier to $\sigma=4.945\times10^{-6}$ moves the nominal control closer to the constrained fuel optimum but produces saturated arcs in $u_y$ and $u_z$. Near an active control boundary, the logarithmic-barrier curvature increases and attenuates the corresponding feedback channel. This prevents aggressive commands through an already saturated actuator but leaves little authority for disturbance rejection. At $\sigma=0.0401$, the controls remain interior and the feedback gains retain sufficient authority to correct the noisy trajectory while respecting the prescribed bounds. This behavior is consistent with the constraint-aware feedback analysis in Ref.~\cite{abhijeet2026safeoptimalcontrolusing}.

\section{Conclusions and Future Work}
\label{sec:conclusions}

This paper presented Box-iLQR, a logarithmic-barrier extension of iLQR for nonlinear optimal-control problems with bounded states and controls. The barrier formulation preserves feasibility, adds positive curvature to the local quadratic model, and is combined with an adaptive continuation strategy to approach the constrained optimum. At a finite barrier value, the backward pass also provides a constraint-aware time-varying feedback policy for closed-loop execution.

Box-iLQR was demonstrated on spacecraft attitude control, a minimum-fuel planar orbit transfer, and an $L_2$-to-$L_2$ Halo-orbit transfer. Across these examples, the method produced feasible constrained trajectories, recovered the expected saturation and bang--off--bang control structures, and achieved high terminal accuracy. Closed-loop simulations further showed that the finite-barrier feedback gains substantially reduce trajectory and terminal errors under disturbances, measurement noise, and control uncertainty. Thus, Box-iLQR provides constrained trajectory optimization and local feedback stabilization within a single framework.

Future work will consider inexact barrier continuation to reduce computational cost, direct treatment of terminal equality constraints, and extension to free-final-time optimal-control problems.

\bibliographystyle{AAS_publication}   
\bibliography{references}   

@misc{abhijeet2026safeoptimalcontrolusing,
      title={Safe Optimal Control using Log Barrier Constrained iLQR}, 
      author={Abhijeet and Suman Chakravorty},
      year={2026},
      eprint={2602.05046},
      archivePrefix={arXiv},
      primaryClass={math.OC},
      url={https://arxiv.org/abs/2602.05046}, 
}

@article{betts1998survey,
  author  = {Betts, John T.},
  title   = {Survey of Numerical Methods for Trajectory Optimization},
  journal = {Journal of Guidance, Control, and Dynamics},
  year    = {1998},
  volume  = {21},
  number  = {2},
  pages   = {193--207},
  doi     = {10.2514/2.4231}
}

@article{rao2009survey,
  author  = {Rao, Anil V.},
  title   = {A Survey of Numerical Methods for Optimal Control},
  journal = {Advances in the Astronautical Sciences},
  year    = {2009},
  volume  = {135},
  number  = {1},
  pages   = {497--528},
  note    = {AAS Paper 09-334}
}

@article{kelly2017introduction,
  author  = {Kelly, Matthew},
  title   = {An Introduction to Trajectory Optimization: How to Do Your Own
             Direct Collocation},
  journal = {SIAM Review},
  year    = {2017},
  volume  = {59},
  number  = {4},
  pages   = {849--904},
  doi     = {10.1137/16M1062569}
}

@article{malyuta2021advances,
  author  = {Malyuta, Danylo and Yu, Yue and Elango, Purnanand and
             A{\c c}{\i}kme{\c s}e, Beh{\c c}et},
  title   = {Advances in Trajectory Optimization for Space Vehicle Control},
  journal = {Annual Reviews in Control},
  year    = {2021},
  volume  = {52},
  pages   = {282--315},
  doi     = {10.1016/j.arcontrol.2021.04.013}
}

@article{vadali1984optimalopenloop,
  author  = {Vadali, Srinivas R. and Junkins, John L.},
  title   = {Optimal Open-Loop and Stable Feedback Control of Rigid
             Spacecraft Attitude Maneuvers},
  journal = {The Journal of the Astronautical Sciences},
  year    = {1984},
  volume  = {32},
  number  = {2},
  pages   = {105--122}
}

@article{carrington1986optimal,
  author  = {Carrington, C. K. and Junkins, John L.},
  title   = {Optimal Nonlinear Feedback Control for Spacecraft Attitude
             Maneuvers},
  journal = {Journal of Guidance, Control, and Dynamics},
  year    = {1986},
  volume  = {9},
  number  = {1},
  pages   = {99--107},
  doi     = {10.2514/3.20073}
}

@article{phogat2018discretetime,
  author  = {Phogat, Karmvir Singh and Chatterjee, Debasish and
             Banavar, Ravi N.},
  title   = {Discrete-Time Optimal Attitude Control of a Spacecraft with
             Momentum and Control Constraints},
  journal = {Journal of Guidance, Control, and Dynamics},
  year    = {2018},
  volume  = {41},
  number  = {1},
  pages   = {199--211},
  doi     = {10.2514/1.G002861}
}

@article{haberkorn2004lowthrust,
  author  = {Haberkorn, Thomas and Martinon, Pierre and Gergaud, Joseph},
  title   = {Low-Thrust Minimum-Fuel Orbital Transfer: A Homotopic Approach},
  journal = {Journal of Guidance, Control, and Dynamics},
  year    = {2004},
  volume  = {27},
  number  = {6},
  pages   = {1046--1060},
  doi     = {10.2514/1.4022}
}

@article{ross2007lowthrust,
  author  = {Ross, I. Michael and Gong, Qi and Sekhavat, Pooya},
  title   = {Low-Thrust, High-Accuracy Trajectory Optimization},
  journal = {Journal of Guidance, Control, and Dynamics},
  year    = {2007},
  volume  = {30},
  number  = {4},
  pages   = {921--933},
  doi     = {10.2514/1.23181}
}

@article{taheri2016enhancedsmoothing,
  author  = {Taheri, Ehsan and Kolmanovsky, Ilya and Atkins, Ella M.},
  title   = {Enhanced Smoothing Technique for Indirect Optimization of
             Minimum-Fuel Low-Thrust Trajectories},
  journal = {Journal of Guidance, Control, and Dynamics},
  year    = {2016},
  volume  = {39},
  number  = {11},
  pages   = {2500--2511},
  doi     = {10.2514/1.G000379}
}

@article{mayne1966secondorder,
  author  = {Mayne, David Q.},
  title   = {A Second-Order Gradient Method for Determining Optimal
             Trajectories of Non-Linear Discrete-Time Systems},
  journal = {International Journal of Control},
  year    = {1966},
  volume  = {3},
  number  = {1},
  pages   = {85--95},
  doi     = {10.1080/00207176608921369}
}

@inproceedings{li2004iterative,
  author    = {Li, Weiwei and Todorov, Emanuel},
  title     = {Iterative Linear Quadratic Regulator Design for Nonlinear
               Biological Movement Systems},
  booktitle = {Proceedings of the First International Conference on
               Informatics in Control, Automation and Robotics},
  year      = {2004},
  volume    = {1},
  pages     = {222--229},
  publisher = {SciTePress},
  doi       = {10.5220/0001143902220229}
}

@inproceedings{tassa2012synthesis,
  author    = {Tassa, Yuval and Erez, Tom and Todorov, Emanuel},
  title     = {Synthesis and Stabilization of Complex Behaviors through
               Online Trajectory Optimization},
  booktitle = {2012 IEEE/RSJ International Conference on Intelligent Robots
               and Systems},
  year      = {2012},
  pages     = {4906--4913},
  publisher = {IEEE},
  doi       = {10.1109/IROS.2012.6386025}
}

@inproceedings{tassa2014controllimited,
  author    = {Tassa, Yuval and Mansard, Nicolas and Todorov, Emanuel},
  title     = {Control-Limited Differential Dynamic Programming},
  booktitle = {2014 IEEE International Conference on Robotics and Automation},
  year      = {2014},
  pages     = {1168--1175},
  publisher = {IEEE},
  doi       = {10.1109/ICRA.2014.6907001}
}

@book{nocedal2006numerical,
  author    = {Nocedal, Jorge and Wright, Stephen J.},
  title     = {Numerical Optimization},
  edition   = {2},
  publisher = {Springer},
  address   = {New York},
  year      = {2006},
  doi       = {10.1007/978-0-387-40065-5}
}

@article{Saloglu2024AccelerationBased,
  author  = {Saloglu, Keziban and Taheri, Ehsan},
  title   = {Acceleration-Based Switching Surfaces for Impulsive Trajectory Design Between Cislunar Libration Point Orbits},
  journal = {The Journal of the Astronautical Sciences},
  year    = {2024},
  volume  = {71},
  number  = {2},
  pages   = {13},
  doi     = {10.1007/s40295-024-00432-z},
  url     = {https://doi.org/10.1007/s40295-024-00432-z}
}

@INPROCEEDINGS{Yamamoto2026NonlinearMPC,
  author={Yamamoto, Koya and Taheri, Ehsan and Junkins, John},
  booktitle={2026 American Control Conference (ACC)}, 
  title={Nonlinear MPC for Cislunar Trajectory Correction Under Impulsive Maneuver Uncertainties: A Computational Performance Analysis}, 
  year={2026},
  volume={},
  number={},
  pages={1449-1454},
  doi={}}

@misc{Yamamoto2025JacobianColoring,
  author       = {Yamamoto, Koya and Taheri, Ehsan and Junkins, John L.},
  title        = {Jacobian Coloring for Multi-Impulse Cislunar Trajectory Optimization},
  year         = {2025},
  howpublished = {SSRN},
  note         = {Available at SSRN: \url{https://ssrn.com/abstract=5545742}},
  doi          = {10.2139/ssrn.5545742}
}

@article{Lantoine2012HybridDDP,
  author  = {Lantoine, Gregory and Russell, Ryan P.},
  title   = {A Hybrid Differential Dynamic Programming Algorithm for Constrained Optimal Control Problems. Part 1: Theory},
  journal = {Journal of Optimization Theory and Applications},
  volume  = {154},
  number  = {2},
  pages   = {382--417},
  year    = {2012},
  doi     = {10.1007/s10957-012-0039-0}
}

@article{Lantoine2012HybridDDPApplication,
  author  = {Lantoine, Gregory and Russell, Ryan P.},
  title   = {A Hybrid Differential Dynamic Programming Algorithm for Constrained Optimal Control Problems. Part 2: Application},
  journal = {Journal of Optimization Theory and Applications},
  volume  = {154},
  number  = {2},
  pages   = {418--442},
  year    = {2012},
  doi     = {10.1007/s10957-012-0038-1}
}

@misc{abhi_DDP_ILQR,
  title        = {A Sequential Quadratic Programming Perspective on Optimal Control},
  author       = {Abhijeet and Chakravorty, Suman},
  year         = {2026},
  eprint       = {2605.25318},
  archivePrefix= {arXiv},
  primaryClass = {math.OC},
  url          = {https://arxiv.org/abs/2605.25318}
}

@INPROCEEDINGS{abhi_DDP_ILQR_ACC,
  author={Abhijeet and Chakravorty, Suman},
  booktitle={2026 American Control Conference (ACC)}, 
  title={A Sequential Quadratic Programming Perspective on Optimal Control}, 
  year={2026},
  volume={},
  number={},
  pages={3862-3867},
  doi={}}

@inproceedings{Howell2019ALTRO,
  author    = {Howell, Taylor A. and Jackson, Brian E. and Manchester, Zachary},
  title     = {ALTRO: A Fast Solver for Constrained Trajectory Optimization},
  booktitle = {2019 IEEE/RSJ International Conference on Intelligent Robots and Systems (IROS)},
  pages     = {7674--7679},
  year      = {2019},
  publisher = {IEEE},
  doi       = {10.1109/IROS40897.2019.8967788}
}

@inproceedings{Vanroye2023FATROP,
  author    = {Vanroye, Lander and Sathya, Ajay and De Schutter, Joris and Decr{\'e}, Wilm},
  title     = {FATROP: A Fast Constrained Optimal Control Problem Solver for Robot Trajectory Optimization and Control},
  booktitle = {2023 IEEE/RSJ International Conference on Intelligent Robots and Systems (IROS)},
  pages     = {10036--10043},
  year      = {2023},
  publisher = {IEEE},
  doi       = {10.1109/IROS55552.2023.10342336}
}

@article{Wachter2006IPOPT,
  author  = {W{\"a}chter, Andreas and Biegler, Lorenz T.},
  title   = {On the Implementation of an Interior-Point Filter Line-Search Algorithm for Large-Scale Nonlinear Programming},
  journal = {Mathematical Programming},
  volume  = {106},
  number  = {1},
  pages   = {25--57},
  year    = {2006},
  doi     = {10.1007/s10107-004-0559-y}
}

\end{document}